\documentclass[12pt]{article}
\usepackage[utf8]{inputenc}
\usepackage{comment}
\usepackage{algorithm}
\usepackage{algorithmic}
\usepackage{mdwlist}
\usepackage{amsfonts}
\usepackage{ulem}
\usepackage{amsmath} 
\usepackage{amsxtra}
\usepackage{mathtools}
\usepackage{fancyhdr}

\newtheorem{theorem}{Theorem}
\usepackage{titlesec}
\usepackage{indentfirst}
\usepackage{booktabs}
\usepackage{verbatim}
\usepackage{color}
\usepackage{latexsym}
\usepackage{comment}
\usepackage{amscd}
\usepackage{esvect}
\usepackage{url}
\usepackage{upgreek}
\usepackage{caption}   
\usepackage[page,header]{appendix}
\usepackage{mathrsfs}
\usepackage{float}
\usepackage{lipsum}
\usepackage[export]{adjustbox} 
\usepackage{subfiles}
\usepackage{subfigure}
\usepackage{bm}
\usepackage{geometry}
\usepackage{epstopdf}
\usepackage{colortbl} 
\usepackage{color}    
\usepackage[cal=cm,scr=rsfs]{mathalfa}

\usepackage[x11names]{xcolor}

\begin{document}

\title{Efficient stochastic asymptotic-preserving scheme for tumor growth models with uncertain parameters}

\date{}
\maketitle

\begin{center}
    \author{
    Ning Jiang $^{\dagger}$,\quad 
    Liu Liu $^{\ddagger}$,\quad 
    Huimin Yu$^{\S}$}
\end{center}

\footnote{\textsuperscript{$\dagger$} School of Mathematics and Statistics, Wuhan University, Wuhan 430072, P. R. China (njiang@whu.edu.cn).}

\footnote{\textsuperscript{$\ddagger$} The Chinese University of Hong Kong, Hong Kong (lliu@math.cuhk.edu.hk)}

\footnote{\textsuperscript{$\S$} School of Mathematics and Statistics, Wuhan University, Wuhan 430072, P. R. China (yuhuimin@whu\\.edu.cn).}

\footnote{Funding: N. Jiang acknowledges the support by NSFC grants 12371224, 11971360, 11731008 and the Strategic Priority Research Program of Chinese Academy of Sciences grant XDA25010404. L. Liu acknowledges the support by National Key R\&D Program of China (2021YFA1001200), Ministry of Science and Technology in China, Early Career Scheme (24301021) and General Research Fund (14303022 \& 14301423) funded by Research Grants Council of Hong Kong.}

\begin{abstract}
   In this paper, we investigate a class of tumor growth models governed by porous medium-type equations with uncertainties arisen from the growth function, initial condition, tumor support radius or other parameters in the model. We develop a stochastic asymptotic preservation (s-AP) scheme in the generalized polynomial chaos-stochastic Galerkin (gPC-SG) framework, which remains robust for all index parameters $m\geq 2$. The regularity of the solution to porous medium equations in the random space is studied, and we show the s-AP property, ensuring the convergence of SG system on the continuous level to that of Hele-Shaw dynamics as $m \to \infty$. Our numerical experiments, including capturing the behaviors such as finger-like projection, proliferating, quiescent and dead cell's evolution, validate the accuracy and efficiency of our designed scheme. The numerical results can describe the impact of stochastic parameters on tumor interface evolutions and pattern formations.



\end{abstract}

\section{Introduction}

Mathematical modeling has been an essential tool in cancer research, providing valuable insights on tumor growth and progression. Many existing work focused on model formulation and analysis to predict tumor development. 
Various models have been proposed to describe the tumor's behaviour, including stochastic models based on reaction-diffusion equations \cite{greenspan1976growth}, phase field models based on Cahn–Hilliard equations \cite{garcke2016cahn}, and mechanical models using porous medium equations \cite{perthame2014hele}. For further details, readers can refer to textbooks \cite{cristini2017introduction,cristini2011multiscale} and review articles \cite{araujo2004history, byrne2006modelling, lowengrub2009nonlinear, perthame2014hele}.

Tumor growth is a highly complex biological process, evolving through distinguishable phases and being affected by numerous factors. Many mathematicians have been devoted to modeling and analyzing individual and synergistic
effects, such as nutrient concentration \cite{feng2023tumor, jacobs2023tumor} and vascularization \cite{cristini2003nonlinear}. 
The finger-like projections growth patterns in tumor growth have been described by various mathematical models \cite{chen1994spectrum,cristini2010multiscale}. The reaction-diffusion equation model emphasizes the growth dynamics of tumors in a confined space and interactions between cells, revealing the impact of Laplace instability on border tumor morphology and leading to finger-like structures \cite{gatenby1996reaction,greenspan1976growth}. The Cahn-Hilliard equation 
is suitable for simulating the evolution of tumor cells at the microscopic scale, capturing the interface phenomena between tumors and surrounding tissue and can capture the instability and curvature-driven behavior of phase transition interfaces \cite{chen1994spectrum,garcke2016cahn}. 
Nevertheless, multiscale models combine biological processes at different scales, we refer to 
\cite{ cristini2017introduction, cristini2011multiscale} for a review. The porous medium model, in particular, considers both reaction and diffusion processes and is suitable to describe tumor growth in regions with restricted blood flow and low oxygen supply \cite{perthame2014hele, saffman1958penetration}. 

We consider porous medium-type equations in this work. These models are categorized by the physical parameter \( m \) that governs different constitutive relations, which connect pressure $p$ with density $\rho$
through $p(\rho)=\frac{m}{m-1}\rho^{m-1}$ for $m\gg 1$. The nonlinearity and degeneracy in the diffusion bring significant challenges in numerical simulations, not to say capturing the singular free boundary limit. Previous studies have shown that porous medium type-equations exhibit asymptotic behavior to Hele-Shaw dynamic as the parameter \(m\) approaches infinity \cite{aronson1998limit,gil2001convergence,igbida2002mesa,kim2003uniqueness}.
The authors in earlier works \cite{liu2021toward,liu2018accurate} have developed an asymptotic-preserving (AP) numerical scheme based on a prediction-correction reformulation that can accurately approximate front propagation, with the semi-discrete scheme convergent to the free boundary limit equation as $m\to\infty$. 

Studying the uncertainty quantification (UQ) problems, on the other hand, is also important to understanding the complicated biological behaviour of tumor growth. Popular approaches to solve forward UQ problems include Monte Carlo simulation \cite{zhang2021modern}, generalized polynomial chaos (gPC)-based  stochastic Galerkin (SG) and stochastic collocation (SC) methods \cite{ghanem2003stochastic, lemaitre2010spectral, Xiu, Xiu2002}. In addition, to validate and predict the model given some observation data, Bayesian inference method and other deep learning methods have been studied, we mention some recent work \cite{Falcó2023, Feng2024, Kahle2019, Kostelich2011, Paixao2021, Selvanambi2020,Zhang2019}. 

\textbf{Main contributions.} In this work, we study a class of porous medium-type tumor growth models that contain stochastic parameters arisen from growth function, initial data or radius support of tumors. We design initial values with finger-like properties to investigate the propagation of uncertainties on interface evolution, pattern changes and tumor growth, providing some new perspectives on modeling interfacial instability and tumor growth dynamics that contain random uncertainties. An efficient numerical method that satisfies the stochastic asymptotic-preservation (s-AP) property \cite{jin2015sAP} is proposed, based on the AP scheme developed in \cite{liu2018accurate} for deterministic problem, and we adopt framework of the gPC-based stochastic Galerkin (SG) method. Regarding analysis, we study the regularity of solution in the random space, and show the s-AP property of the gPC-SG system without discretization in time and space, that is, our semi-discrete SG system converges to the corresponding SG system of the limit equation as $m\to\infty$. 

The paper is organized as follows. In Section 2, we introduce a class of tumor growth models with uncertain parameters, related to vitro and vivo nutrient models. We analyze regularity of the solution in the random space and stochastic-AP property for the SG system without discretization in time and space. In Section 3, we develop a fully discrete scheme in the gPC-SG framework for the porous media equation of vivo nutrient model. In Section 4, a series of numerical examples {that consider both 1D and 2D random variables} are conducted to demonstrate the accuracy and efficiency of our proposed scheme, with problems designed to investigate finger-like invasion characteristics and various cell models. 

\section{Tumor growth model}

 We are interested in a class of mechanical tumor growth models, specifically those governed by porous medium type-equations indexed by a physical parameter \(m\), which characterizes the relationship between the pressure and the density \cite{benilan1996singular}. As \( m \to \infty \), these equations asymptotically approach the natural Hele-Shaw dynamics. 

There may exist uncertainties coming from growth function \(g\), initial data\(f\), radius of the support of \(\rho\) or other model parameters. The random variable \(z\) is a \( n \)-dimensional vector with support \( I_z\) characterizing the random uncertainties in the system. We assume that it has a prescribed probability density function \( \pi(z) \ge 0 \). 

Let \( \Omega \) be a bounded open set in \( \mathbb{R}^2 \), representing the domain of tumor growth. For \( T > 0 \), define \( Q_T := \Omega \times (0, T) \) and \( \Sigma_T := \partial \Omega \times (0, T) \). The continuity of mass reads
 \begin{equation} \label{eq:main}
 	\frac{\partial \rho}{\partial t} + \nabla \cdot (\rho u) = g(x, t, \rho,z).
 \end{equation}
 
 Denote by \( \rho(x, t, z) \) the cell population density, which is transported by a velocity field \( u(x,t,z) \) and influenced by a growth rate function \( g(x, t,\rho, z) \). The velocity \( u(x,t,z) \) is assumed to follow Darcy's law, \( u = -\nabla p \), where the pressure \( p \) satisfies the power law \(p(x,t,z) = \frac{m}{m-1} \rho(x,t,z)^{m-1}\), for $m \geq 2$. The tumor boundary expands with a finite normal speed \( s = -\nabla p \cdot \mathbf{n}|_{\partial D} \), where \( \mathbf{n} \) is the outer normal vector on the boundary. In this paper, we consider the growth function in the form:

\begin{equation} \label{eq:growth_function}
	g(x, t, \rho,z) = h(x,c,z) \rho, \quad 0 < h(x,c,z) \in L^\infty(\Omega),
\end{equation}
where \( h(x,c,z) \) represents the growth rate function and $c(x,z)$ is nutrition density in vitro or vivo equation \cite{liu2019analysis,perthame2014hele}, reflecting the tumor micro-environment. In addition, since we focus on early-stage tumor development, cell apoptosis is negligible, thus \( h(x,c,z) \) is  strictly positive. We assume that \( \rho \) vanish on \( \Sigma_T \). Let \( f(x,z) \) denote the initial condition, which is typically an arbitrary function taking values in \( [0, 1] \). For any \( m \geq 2 \), the evolution of the tumor density satisfies the following system:

\begin{equation} \label{eq:tumor_system}
	(P_m) \quad 
	\begin{cases} 
	\displaystyle	\frac{\partial \rho}{\partial t} -\frac{m}{m-1}  \nabla \cdot \left \{ {\rho \nabla(\rho^{m-1 })} \right \}  = h(x, c,z)\rho, & \text{on } Q_T, \\ 
		\rho(x,t,z) = 0 & \text{on } \Sigma_T, \\ 
		\rho(x, 0,z) = f(x,z) & \text{on } \Omega.
	\end{cases}
\end{equation}

\subsection{Regularity in z }

In this section, we show that the solution $\rho$ of \eqref{eq:tumor_system} preserves the regularity of the initial function in the random space. The stability of the solution is obtained, in a suitable weighted Sobolev norm. For each $z$, 
we follow the same regularity assumptions as in \cite[Section 3.1]{feng2024regularity}. We summarize the result in the following Theorem. 
     
\begin{theorem}
    Assume that for $m >2$, initial data 
    $\|f_0\|_{W^{{s}}(0)} \leq  \beta$, under the same assumptions as in \cite{feng2024regularity} for each $z$, we define the $W^s(t)$ norm of a function $f$ as: 
	\[
    {
		\|f\|^2_{W^s(t)} = \sum_{l=0}^s  \sum_{k=0}^{l} \binom{l}{k} \|\partial_z^{l - k} f\|^2_{H^1}, 
	}
	\]
    {where \( \binom{l}{k} \) denotes the binomial coefficient, with \( 0 \leq k \leq l \) and \( 0 \leq l \leq s \).} Then we have
    \[
\|\rho \|_{W^{{s}}(t)} \leq \beta e^{\tilde{C}t/2},
\]
here $\tilde{C}$ is a constant.
\end{theorem}

\textit{proof}: The expression for \( p \) allows the flux \( -\frac{m}{m-1}  \nabla \cdot \left \{ {\rho \nabla(\rho^{m-1 })} \right \}\) to be equivalently written as \( -\Delta(\rho^m)  \). We first take \( l \)-th derivative with respect to \( z \) to (\ref{eq:tumor_system}):
	\[
	\frac{\partial^l}{\partial z^l} \left( \frac{\partial \rho}{\partial t} \right) = \Delta \left( m \rho^{m-1} \frac{\partial^l}{\partial z^l} \rho \right) + {\sum_{k=0}^l \binom{l}{k} \frac{\partial^k}{\partial z^k} h(x, c, z) \cdot \frac{\partial^{l-k}}{\partial z^{l-k}} \rho}.
	\]
 Multiplying both sides by \( \frac{\partial^l}{\partial z^l} \rho \) and integrating on $Q_T\times I_z$, one gets
 \begin{align*}
& \int_{Q_T \times I_z} \frac{\partial^l}{\partial z^l} \left( \frac{\partial \rho}{\partial t} \right) \cdot \frac{\partial^l}{\partial z^l} \rho \, dx \, dy \, dz 
=  \int_{Q_T \times I_z} \Delta \left( m \rho^{m-1} \frac{\partial^l}{\partial z^l} \rho \right) \cdot \frac{\partial^l}{\partial z^l} \rho \, dx \, dy \, dz, \\
& + {\sum_{k=0}^l \binom{l}{k} \int_{Q_T \times I_z} \frac{\partial^k}{\partial z^k} h(x, c, z) \cdot \frac{\partial^{l-k}}{\partial z^{l-k}} \rho \cdot \frac{\partial^l}{\partial z^l} \rho \, dx\, dy\, dz}.
\end{align*}
Integral by part and sum over {$l$}, then right-hand side becomes:
		\begin{align*}
			&\sum_{l=0}^{{s}} \int_{Q_T \times I_z} \Delta \left( m \rho^{m-1} \right) \cdot \left( \frac{\partial^l}{\partial z^l} \rho \right)^2 \, dx \, dy \, dz \\
			& + \sum_{l=0}^{{s}} \sum_{k=0}^l \binom{l}{k} \int_{Q_T \times I_z} \frac{\partial^k}{\partial z^k} h(x, c,z) \cdot \frac{\partial^{l-k}}{\partial z^{l-k}} \rho \cdot \frac{\partial^l}{\partial z^l} \rho \, dx \, dy \, dz \\
			&\leq C \sum_{l=0}^{{s}} \left\| \frac{\partial^l}{\partial z^l} \rho \right\|_{{H^1}}^2 + {C} \sum_{l=0}^{{s}} \sum_{k=0}^l \binom{l}{k} \left\| \frac{\partial^k}{\partial z^k} h \right\|_\infty \left\| \frac{\partial^{l-k}}{\partial z^{l-k}} \rho \right\|_{L^2} \left\| \frac{\partial^l}{\partial z^l} \rho \right\|_{L^2},\\
			&\leq C {\sum_{l=0}^{s}} \left\| {\frac{\partial^l}{\partial z^l} \rho} \right\|_{H^1}^{2} + \frac{{C}}{2} {\sum_{l=0}^{s}} \sum_{k=0}^l \binom{l}{k} \left\| \frac{\partial^k}{\partial z^k} h \right\|_\infty \left( \left\| \frac{\partial^{l-k}}{\partial z^{l-k}} \rho \right\|_{{H^1}}^2 + \left\| \frac{\partial^l}{\partial z^l} \rho \right\|_{{H^1}}^2 \right).
		\end{align*}
        
\noindent Combine the summation on the right-hand-side, we have
\[
\partial_t \|\rho \|^2_{W^{{s}}(t)} \leq \tilde{C} \|\rho \|^2_{W^{{s}}(t)}, 
\]
where $C$ and $\tilde{C}$ are constants. 
Applying the Gr\"onwall’s inequality, one has 
\[
\|\rho \|_{W^{{s}}(t)} \leq e^{\tilde{C}t/2} \|f _0\|_{W^{{s}}(0)} \leq \beta e^{\tilde{C}t/2}.
\]

\subsection{The gPC-SG method}
\label{sec:UQ}

Among various numerical methods for solving UQ problems, the generalized polynomial chaos (gPC)-based stochastic Galerkin (SG) method have been used popularly and shown successful in broad applications \cite{Xiu}. It is computationally efficient and can achieve spectral accuracy in the random space, provided the solution is smooth enough with respect to the random variable. One inserts the solution ansatz $\rho_K$ defined by:
\begin{equation}
\rho(x, t,z) \approx \rho_K(x, t,z) =
\sum_{|{k}|=1}^{K} \rho_{{k}}( x,t) \psi_{{k}}(z) = \boldsymbol{\rho} \cdot \boldsymbol{\psi}, \quad K =
\binom{n + P}{n}, 
\end{equation}
where \(n\) is the dimension of the random variable \(z\), \(P\) is polynomial order, \( {k} = (k_1, \dots, k_n) \) is a multi-index for the vector with \( |{k}| = k_1 + \cdots + k_n \), with
\begin{equation}\label{gPC-vec}\boldsymbol{\rho} = (\rho_1, \dots, \rho_K),\quad \boldsymbol\psi = (\psi_1, \dots, \psi_K). \end{equation}
Here \(\{ \boldsymbol{\psi}{_{{k}}(z)} \} \) are the orthonormal basis functions that form \( \mathbb{P}_P^n \)  
(the set of \( n \)-variate orthonormal polynomials of degree up to \( P \geq 1 \)) and satisfy  
\[
\int_{I_z} \psi_{{k}}(z) \psi_{{l}}(z) \pi(z) \, dz = \delta_{kl}, \quad 1 \leq |{k}|,\, |{l}| \leq K = \dim(\mathbb{P}_P^n),
\]

We introduce some notations for the space, inner produce and norm that will be used:  
$$
H = L^2(I_z; \pi(z) \,dz),\quad 
\langle f, g \rangle_H = \int_{I_z} f g \,\pi(z) \,dz,\quad 
\| f \|_H = \left( \int_{I_z} f^2 \pi(z) \,dz \right)^{\frac{1}{2}}.
$$

Apply \(\rho_K\) to (\ref{eq:tumor_system}) and perform a standard Galerkin projection, we get
\begin{equation} 
	\partial_t \rho_i + \nabla \cdot (\rho_i  u_i)   = \sum_{j = 1}^{K}  H_{ij}
	\rho_i, 
\end{equation}
where the matrix $\mathbf{H}(x,c) := (H_{ij})_{K\times K}$ is defined by 
\begin{equation}
\label{HHH}
	H_{ij}(x,c) = \int_{I_z} h(x,c, z) \psi_i(z) \psi_j(z) \, dz.
\end{equation}
The initial values of each component of \({\rho(x,t,z)}\) are given by 
\[\rho_k(x,0) = \langle f(x), \psi_k \rangle_H, \quad |k| = 1, \dots, K.\]

We now write it in a vector form: 
\begin{equation}
\label{eq:vectorized_system}
		\begin{cases}
			\displaystyle \frac{\partial \boldsymbol{\rho}}{\partial t}  + \nabla \cdot (\boldsymbol{\rho}  \circ \mathbf{u}) = \mathbf{H} \boldsymbol{\rho} & \text{in } (0, T) \times \mathcal{X}, \\
			\boldsymbol{\rho} = 0 & \text{on } (0, T) \times \partial\mathcal{X}, \\
			\boldsymbol{\rho}(x,0) = \mathbf{f}(x) & \text{in } \mathcal{X}.
		\end{cases}
\end{equation}

The Hadamard (elementwise) product is denoted by \( \circ \), with \(\boldsymbol{\rho}(t,x) = (\rho_1, \dots, \)

\(\rho_K)^T,  \boldsymbol{u} = \begin{pmatrix}
		u_1, \dots, u_K
	\end{pmatrix}^\top, \boldsymbol{f} = \begin{pmatrix}
	f_1, \dots, f_K
	\end{pmatrix}^\top \in \mathcal{X} = (L^1(\Omega))^{K}\).

\subsection{Stochastic-AP property}

For simplicity, we denote the flux term as \( \Delta(\rho^m)  \) in the model. The SG system (\ref{eq:vectorized_system}) is given by 
\begin{equation}
    \label{eq:system}
    \begin{cases}
    \displaystyle
        \frac{\partial \boldsymbol{\rho}}{\partial t}  = \mathbf{\Delta} (\mathbf{\rho}^{m}) + \mathbf{H} \mathbf{\rho}, & \text{in } (0, T) \times \mathcal{X}, \\[2mm]
        \mathbf{\rho} = 0, & \text{on } (0, T) \times \partial \mathcal{X}, \\[2mm]
        \mathbf{\rho}(x,0) = \mathbf{f}(x), & \text{in } \mathcal{X},
    \end{cases}
\end{equation}
	where \(\mathbf{\Delta}\) acts componentwise:  	\((\mathbf{\Delta} \boldsymbol{\rho})_i = \Delta \rho_i\) for each \(i\), 
and $\overset{\sim}{\mathcal{D}} = \{x \in \mathcal{X} \mid \hat{\mathbf{f}}(x) \neq \mathbf{1}\}$, 
$\hat{\mathcal{D}} = \{x \in \mathcal{X} \mid \mathbf{f}(x) \neq \mathbf{1}\}$. 
 
     Based on \cite[Lemma 2]{benilan1996singular}, as \( m\to\infty \), the limiting system for  (\ref{eq:vectorized_system}) is written as
    
	\begin{equation}\label{eq:vectorized_limit}
		\begin{cases}
        \displaystyle
			\frac{\partial \boldsymbol{\rho}}{\partial t}  = \mathbf{\Delta} \mathbf{w} + \mathbf{H}\boldsymbol{\rho} & \text{in } (0, T) \times \mathcal{X}, \\[2mm]
			0 \leq \boldsymbol{\rho} \leq \mathbf{1}, \quad \mathbf{w} \geq 0, &\text{in } (0, T) \times \mathcal{X}, \quad\\[2mm] 
			\mathbf{w} = 0 & \text{on } \partial\mathcal{X}  \quad\text{and in} \quad(0, T) \times (\overset{\sim}{\mathcal{D}}), \\[2mm]
			\boldsymbol{\rho}(x,0) = \mathbf{\hat{f}}(x) & \text{in } \mathcal{X},
		\end{cases}
	\end{equation}
where \(\mathbf{\hat{f}(x)} = \mathbf{f}(x) \chi_{[\mathbf{\underline{w}}=\mathbf{0}]} + \chi_{[\mathbf{\underline{w}}>\mathbf{0}]}
\), \(\mathbf{\underline{w}}(t,x) = (\underline{w}_1, \dots, \underline{w}_K)^T\) satisfies the SG mesa problem:
	\begin{equation*}
		\mathbf{\underline{w}} \in H^1_0(\Omega)^{K}, \quad \mathbf{{\Delta} \underline{w}} \in L^1(\Omega)^{K}, \quad 0 \leq \mathbf{{\Delta} \underline{w}} + \mathbf{f} \leq \mathbf{1},
	\end{equation*}
	\begin{equation*}
		\mathbf{{\Delta} \underline{w}} + \mathbf{f} - \mathbf{1} = 0\quad\text{a.e. in } (0, T) \times (\hat{\mathcal{D}}), \quad \mathbf{\underline{w}} = 0 \quad\text{a.e. in } (0, T) \times (\hat{\mathcal{D}}\cap\mathcal{X}).
	\end{equation*}

\begin{theorem}
    Under the same assumptions as Theorem 1, as \( m \to \infty \), let \( \boldsymbol{\rho}_m \), \( \boldsymbol{\rho}_\infty \) be the solution to the SG system (\ref{eq:system}) and (\ref{eq:vectorized_limit}), then 
	\begin{enumerate}
		\item \( \boldsymbol{\rho}_m \to \boldsymbol{\rho}_\infty \) in \( C((0, T); L^1(\Omega))^{K})\) as \( m \to \infty \).
		\item Assuming when \(\boldsymbol{\rho} = 1\), \( \mathbf{H}\boldsymbol{\rho}  \leq \tilde{h} \) in \( (L^1(\mathcal{X})) \), with \( \tilde{h} \in L^2_{\text{loc}}((0,T), H^{-1}(\mathcal{X}))\), there exists a unique pair \( (\boldsymbol{\rho}_{\infty}, \mathbf{w}) \) solving the limit system \( (P_\infty) \), where
		\begin{equation}\label{eq:limit_system}
			(P_\infty) \begin{cases}
				\boldsymbol{\rho}_{\infty} \in C([0,T); (L^1(\Omega)^{K})), \quad \mathbf{w} \in L^2_{\text{loc}}((0,T), H_0^1(\mathcal{X})), \\[2mm]
				\boldsymbol{\rho}_{\infty}(x,0) = \mathbf{f}(x), \quad 0 \leq \boldsymbol{\rho} \leq \mathbf{1}, \quad \mathbf{w} \geq 0, \\[2mm]
				\mathbf{w} = 0 \quad \text{in} \quad L^2_{\text{loc}}((0,T), H_0^1((\hat{\mathcal{D}}\cap\mathcal{X})), \\[2mm]
                \displaystyle
				\frac{\partial \boldsymbol{\rho}_{\infty}}{\partial t} = \boldsymbol{\Delta} \mathbf{w} + \mathbf{H} \boldsymbol{\rho}_{\infty} \quad \text{in } (0,T)\times L^1(\mathcal{X}).
			\end{cases}
		\end{equation}
	\end{enumerate}
\end{theorem}

We observe that the only difference between the gPC-SG system for porous medium type equation and its limiting Hele-Shaw dynamics is the flux term \(\mathbf{\Delta} (\boldsymbol{\rho}^{m})\) and \(\mathbf{\Delta} \mathbf{w}\). Similar to the proof in \cite[Theorem 2]{benilan1996singular}, one can easily show the limiting behaviour for the flux term, see Appendix. As \(m \to \infty\), the SG system for the uncertain porous medium model (\ref{eq:system}) satisfies the stochastic-AP property, in particular, it automatically becomes a SG approximation for the limiting stochastic Hele-Shaw dynamics (\ref{eq:vectorized_limit}).

\section{Numerical scheme}

\subsection{2D Radial Symmetric Model with Vitro and Vivo Model}

We consider two specific nutrition models: the vitro  and vivo models studied in \cite{liu2019analysis,perthame2014hele}. Uncertain nutrition density $c(x,z)$ is assumed. For the two-dimensional vitro model, the nutrition equation of $c(x,z)$ is given by
\begin{align}
	- \Delta c + \psi(\rho) c &= 0, \quad x \in D(t,z), \label{eq:invitro_main}\\
	c &= c_B, \quad x \in \mathbb{R}^2 \cap {D(t,z)}^c, 
    \label{eq:invitro_bc}
\end{align}
where $D(t,z) = \{x \in 
\mathbb{R}^2 \mid \rho(x, t,z) > 0\}$ and $\psi(\rho) \geq 0$ with $\psi(0) = 0$. For simplicity, we assume $\psi(\rho) = \chi_{D(t,z)}$. With the region $D(t,z)$, radius $R(t,z)$, for $x \in D(t,z) = B_{R(t,z)}(t,z)$ we have
\begin{equation}
	- \frac{1}{r} \partial_r \left( r \partial_r c \right) + c = 0, 
\end{equation}
Here
\begin{equation}
	c(x,z) =c(r,z) =
	\begin{cases}
		\frac{c_B}{I_0(R(t,z))} I_0(r), & r \in [0, R(t,z)], \\
		c_B, & r > R(t,z),
	\end{cases}
\end{equation}
where $I_0(r)$ is the modified Bessel function of the first kind. 

For the vivo model, the nutrition equation is given by 
\begin{equation}
	- \Delta c + \psi(\rho)\chi_{D(t,z)}c = \chi_{\{\rho = 0\}}(c_B - c). \label{eq:invivo_main}
\end{equation}
For $x \in D(t,z)$, we obtain
\begin{equation}
	- \frac{1}{r} \partial_r \left( r \partial_r c \right) + c = 0,
\end{equation}
for $x \in \mathbb{R}^2 \cap {D(t,z)}^c $:
\begin{equation}
	- \frac{1}{r} \partial_r \left( r \partial_r c \right) = c_B - c,
\end{equation}
giving the solution:
\begin{equation}
	c(x,z) =c(r,z) = c_B + a_1(z) K_0(r) + a_2(z) I_0(r),
\end{equation}
where $K_0(r,z)$ is the modified Bessel function of the second kind. The far-field assumption $c \to c_B$ as $r \to \infty$ gives us $a_2 = 0$. Continuity of $c$ and $\partial_r c$ at $r = R(z)$ leads to 
\begin{align*}
	a_0(z) &= \frac{c_B K_1(R(z))}{K_0(R(z)) I_1(R_z) + K_1(R(z)) I_0(R(z))}, \\
	a_1(z) &= - \frac{c_B I_1(R(z))}{K_0(R(z)) I_1(R(z)) + K_1(R(z)) I_0(R(z))}.
\end{align*}
Thus
\begin{equation}
	c(x,z) =c(r,z) =
	\begin{cases}
		\frac{c_B K_1(R(z))}{K_0(R(z)) I_1(R_z) + K_1(R_z) I_0(R(z))} I_0(r), & r \in [0, R(t,z)], \\[6pt]
		c_B - \frac{c_B I_1(R(z))}{K_0(R(z))I_1(R(z)) + K_1(R(z))I_0(R(z))} K_0(r), & r > R(t,z). 
	\end{cases}\label{c(x,z)}
\end{equation}
With the above assumptions, for any \(m \geq 2\), \(G_0\) is a constant and we set \(h(x,c,z) = G_0 c(r,z)\). Apply (\ref{eq:vectorized_system}), we derive that
\begin{equation}
\label{test I and II}
\frac{\partial \boldsymbol{\rho}}{\partial t} + \nabla \cdot (\boldsymbol{\rho}  \circ \mathbf{u}) = G_0 \mathbf{C} \boldsymbol{\rho},
\end{equation}
where \(\boldsymbol{\rho} = \begin{pmatrix}
		\rho_1, \cdots, \rho_K
	\end{pmatrix}^\top\), \( \boldsymbol{u} = \begin{pmatrix}
		u_1, \cdots, u_K
	\end{pmatrix}^\top
\). The coefficient matrix $\mathbf{C}(x,z) = (C_{ij})_{K\times K}$ is defined as
\begin{equation}
\label{nutrition caculation}
	C_{ij}(x) = \int_{I_z} \mathbf{C}(x, z) \psi_i(z) \psi_j(z) \, dz.
\end{equation}

\subsection{Proliferating, quiescent and dead cells model}
We introduce a biologically more realistic model in \cite{perthame2010some}. Let $\rho_P(x,t)$, $\rho_Q(x,t)$, and $\rho_D(x,t)$ be the cell densities for proliferating, quiescent and dead cells, respectively. We set \(h(x,c,z) = G_0 c(x,z)\) and the cell evolution is governed by the following equations
\begin{equation}
\label{test III}
	\left\{
	\begin{aligned}
		\frac{\partial \rho_P}{\partial t} + \nabla \cdot(\rho_P u) &= G_0 c (x,z)\rho_P - a \rho_P + b \rho_Q, \\
		\frac{\partial \rho_Q}{\partial t} + \nabla \cdot(\rho_Q u) &= a\rho_P-  b \rho_Q - d \rho_Q, \\
		\frac{\partial \rho_D}{\partial t} + \nabla \cdot (\rho_D u) &= d\rho_Q- \mu \rho_D.
	\end{aligned}
	\right.
\end{equation}
The total density is defined by 
\(\rho= \rho_P + \rho_Q + \rho_D, \label{eq:total_density}\), with the velocity field 
\(\displaystyle u = - \nabla p  = -\frac{m}{m - 1} \frac{\partial \rho^{m-1}}{\partial x}\). In this model, the growth rate function \(h(x,c,z)\) consists of the vivo nutritional equation part and the nutritional constants part, that is, $a$, $b$, $d$ and $\mu$. Similarly, one can derive the following gPC-SG system for proliferating, quiescent and dead cells models: 
\begin{equation}
	\left\{
	\begin{aligned}
		\frac{\partial \boldsymbol{\rho}_P}{\partial t} + \nabla \cdot(\boldsymbol{\rho_P} \circ \boldsymbol{u}) &= G_0\mathbf{C} \boldsymbol{\rho}_P - a\boldsymbol{\rho}_P+b \boldsymbol{\rho}_Q \\
		\frac{\partial \boldsymbol{\rho}_Q}{\partial t} + \nabla\cdot(\boldsymbol{\rho_Q} \circ \boldsymbol{u}) &=  a\boldsymbol{\rho}_P-b \boldsymbol{\rho}_Q, -d\boldsymbol{\rho}_Q,\\
		\frac{\partial \boldsymbol{\rho}_D}{\partial t} + \nabla \cdot(\boldsymbol{\rho_D} \circ \boldsymbol{u}) &= d\boldsymbol{\rho}_Q - \mu \boldsymbol{\rho}_D,
	\end{aligned}
	\right.
	\label{eq:biological_model_cases}
\end{equation}
where
\(\boldsymbol{\rho}_P = \begin{pmatrix}
	\rho_{P,1}, \dots, \rho_{P,K}
\end{pmatrix}^\top, 
\boldsymbol{\rho}_Q = \begin{pmatrix}
	\rho_{Q,1}, \dots, \rho_{Q,K}
\end{pmatrix}^\top, 
\boldsymbol{\rho}_D = \begin{pmatrix}
	\rho_{D,1}, \dots, \rho_{D,K}
\end{pmatrix}^\top.\)

\subsection{A fully discretized scheme}

We adopt the AP scheme in \cite{liu2018accurate} for deterministic problems and design the s-AP scheme by using the gPC-SG method. For two-dimensional problem, velocity has two components \( u = (u, v) \), with \( u\) and \( v \) being the velocities along the \( x \) and \( y \) directions, respectively. 

\begin{equation}
u = -\nabla_x p=- \frac{m}{m - 1} \frac{\partial \rho^{m-1}}{\partial x},\qquad
v =-\nabla_y p = - \frac{m}{m - 1} \frac{\partial \rho^{m-1}}{\partial y}.
\label{eq:velocity}
\end{equation}
For each index \(1\le k \le K\), equation (\ref{eq:vectorized_system}) becomes the following equation for velocity \(u_k = (u_k,v_k)\):

\[
\frac{\partial u_k}{\partial t} = m \frac{\partial}{\partial x} \left[ \rho_k^{m-2}  \left( \frac{\partial}{\partial x}(\rho_k u_k) + \frac{\partial}{\partial y}(\rho_k v_k) - G_0 \sum_{l = 1}^{K} C_{lk} \rho_k \right) \right],
\]

\[
\frac{\partial v_k}{\partial t} = m \frac{\partial}{\partial y} \left[ \rho_k^{m-2}  \left( \frac{\partial}{\partial y}(\rho_k v_k) + \frac{\partial}{\partial x}(\rho_k u_k) - G_0 \sum_{l = 1}^{K} C_{lk} \rho_k \right) \right].
\]
We set the computational domain as \( (x, y) \in [a, b] \times [a, b] \). Let \(\displaystyle \Delta x = \frac{b - a}{N_x} \) be the mesh size, and the grid points be
\[
x_i = a + i \Delta x, \quad x_{i+1/2} = a + \left(i + 1/2\right) \Delta x, \quad i \in \{0, 1, \dots, N_x - 1\}, 
\]
and
\[
y_j = a + j \Delta y, \quad y_{j + 1/2} = a + (j + 1/2) \Delta y, \quad \Delta y = \frac{b - a}{N_y}, \quad j \in \{0, 1, \dots, N_y - 1\}.
\]

For each fixed \(k\), we compute both \( {\rho_k} \) and \( ({u_k}, {v_k}) \) on regular grids and specify the half-grid values if needed. We denote
\[
{\rho}_{i,j,k}(t)=\frac{1}{\Delta x \Delta y} \int_{x_{i-1 / 2}}^{x_{i+1 / 2}} \int_{y_{i-1 / 2}}^{y_{i+1 / 2}} \rho_k(x, y, t) \, d x \, d y, \]
\[
u_{i,j,k}(t) \approx u_k\left(x_{i}, y_{j}, t\right), \quad
v_{i,j,k}(t) \approx v_k\left(x_{i}, y_{j}, t\right), 
\]
the discretization for \( {u_k} \) is given by 
\begin{equation}
\begin{aligned}
	\frac{u_{i, j,k}^{n *}-u_{i, j,k}^{n}}{\Delta t}= & \frac{m}{\Delta y} \Bigg\{\left(\rho_{i, j+1 / 2,k}^{n}\right)^{m-2} \frac{\left(\rho^{n}_{k} u^{n *}_{k}\right)_{i, j+1}-\left(\rho^{n}_{k} u^{n *}_{k}\right)_{i, j}}{\Delta y} \\
    &-\frac{1}{2} \left[G_0 {\textstyle \sum_{l=0}^{K}} C_{lk}^{n+1} (\rho{_k^{n}})^{m-1}\right]_{i, j+1} -\frac{1}{2} \left[G_0 {\textstyle \sum_{l=0}^{K}} C_{lk}^{n+1} (\rho{_k^{n}})^{m-1}\right]_{i, j-1} \\
	&  -\left(\rho_{i, j-1 / 2,k}^{n}\right)^{m-2} \frac{\left(\rho^{n}_{k} u^{n *}_{k}\right)_{i ,j}-\left(\rho^{n}_{k} u^{n *}_{k}\right)_{i, j-1}}{\Delta y} \\
	& +\frac{1}{2} \left(\rho_{i, j+1,k}^{n}\right)^{m-2} \frac{\left(\rho^{n}_{k} v^{n *}_{k}\right)_{i+1, j+1}-\left(\rho^{n}_{k} v^{n *}_{k}\right)_{i-1, j+1}}{2 \Delta x} \\
	&  -\frac{1}{2} \left(\rho_{i, j-1,k}^{n}\right)^{m-2} \frac{\left(\rho^{n}_{k} v^{n *}_{k}\right)_{i+1, j-1}-\left(\rho^{n}_{k} v^{n *}_{k}\right)_{i-1, j-1}}{2 \Delta x} \Bigg\}.
\end{aligned}
\end{equation}

\noindent The discretization for $v_k$ is computed similarly.

The value of $\rho$ at half grids are defined by 
\[
\rho_{i+1/2,j,k} = \frac{1}{2} (\rho_{i,j,k} + \rho_{i+1,j,k}), \quad \rho_{i,j+1/2,k} = \frac{1}{2} (\rho_{i,j,k} + \rho_{i,j+1,k}).
\]

After obtaining the values of $u^{n*}_{k}$ and $v^{n*}_{k}$, one can obtain $\rho^{n+1}_{k}$ through the central scheme: 
\[
\begin{aligned}
\rho_{i,j,k}^{n+1} - \rho_{i,j,k}^n 
&= \frac{1}{\Delta x} \left( (F_1)^n_{i+1/2,j,k} - (F_1)^n_{i-1/2,j,k} \right)\\
&+ \frac{1}{\Delta y} \left( (F_2)^n_{i,j+1/2,k} - (F_2)^n_{i,j-1/2,k} \right) 
= \left(G_0 {\textstyle \sum_{l=0}^{K}} C_{lk} \rho{_k}\right)_{i,j}^{n+1} ,
\end{aligned}
\]
where
\[
(F_1)^n_{i\pm1/2,j,k} = \frac{1}{2} \left( \rho_{Lx,k}^n u^{n*}_{k} + \rho_{Rx,k}^n u^{n*}_{k} - |u^{n*}_{k}| (\rho_{Rx,k}^n - \rho_{Lx,k}^n) \right)_{i\pm1/2,j},
\]
\[
(F_2)^n_{i,j\pm1/2,k} = \frac{1}{2} \left( \rho_{Ly,k}^n v^{n*}_{k} + \rho_{Ry,k}^n v^{n*}_{k} - |v^{n*}_{k}| (\rho_{Ry,k}^n - \rho_{Ly,k}^n) \right)_{i,j\pm1/2}.
\]
Here $\rho_{Lx,k}^n, \rho_{Rx,k}^n, \rho_{Ly,k}^n, \rho_{Ry,k}^n$ are half grid values obtained by linear reconstruction:
\[
(\rho_{Lx,k}^n)_{i+1/2,j} = \rho_{i,j,k}^n + \frac{\Delta x}{2} \left( \frac{\partial \rho_{k}}{\partial x} \right)_{i,j}, (\rho_{Rx,k}^n)_{i+1/2,j} = \rho_{i+1,j,k}^n - \frac{\Delta x}{2} \left( \frac{\partial \rho_{k}}{\partial x} \right)_{i+1,j},
\]
\[
(\rho_{Ly,k}^n)_{i,j+1/2} = \rho_{i,j,k}^n + \frac{\Delta y}{2} \left( \frac{\partial \rho_{k}}{\partial y} \right)_{i,j}, (\rho_{Ry,k}^n)_{i,j+1/2} = \rho_{i,j+1,k}^n - \frac{\Delta y}{2} \left( \frac{\partial \rho_{k}}{\partial y} \right)_{i,j+1}.
\]
The discrete gradient is determined by the following expression:
\[
(\frac{\partial \rho_{k}}{\partial x})_{i,j} = \frac{\rho_{i+1,j,k}^n - \rho_{i,j,k}^n}{\Delta x}, \quad (\frac{\partial \rho_{k}}{\partial y})_{i,j} = \frac{\rho_{i,j+1,k}^n - \rho_{i,j,k}^n}{\Delta y}.
\]
In the correction step, we have
\[
u_{i,j,k}^{n+1} = - \frac{m}{m-1} \left( \left( \rho_{i,j,k}^{n} + \rho_{i+1,j,k}^{n} \right)^{m-1} - \left( \rho_{i,j,k}^{n} + \rho_{i-1,j,k}^{n} \right)^{m-1} \right) \frac{1}{2 \Delta x},
\]
\[
v_{i,j,k}^{n+1} = - \frac{m}{m-1} \left( \left( \rho_{i,j,k}^{n} + \rho_{i,j+1,k}^{n} \right)^{m-1} - \left( \rho_{i,j,k}^{n} + \rho_{i,j-1,k}^{n} \right)^{m-1} \right) \frac{1}{2 \Delta y}.
\]
Finally, we solve the $C^{n+1}(x)$ via (\ref{c(x,z)}) and update the nutrition matrix by 
\[
 C_{lk}^{n+1}(x) = \int_{I_z} c(x, z)^{n+1} \psi_l(z) \psi_k(z) \, dz.\]

\section{Numerical examples}

In our numerical experiments, we consider the tumor growth model of porous medium type equations in 2D. We set the domain $x \in \Omega=[-2.2,2.2]\times [-2.2,2.2]$, and assume the no-flux boundary condition $\rho u = 0$ for $x\in \partial\Omega$. In our numerical tests, assume the random variable $z$ satisfies the uniform distribution on $[-1,1]$, thus the Legendre polynomial basis function is used. 

We consider uncertainties arising from the nutrition function \(c(x,z)\), the initial data\(f(x,z)\), the radius of support of \(\rho\), and the parameter data. 
Denote the gPC coefficients by $\rho_k$ \((1 \le k \le K)\), the mean, variance and standard deviation of $\rho$ are define as:
\begin{equation}
	\mathbb{E}[\rho] \approx \rho_1, \quad 
    \mathbb{SD}[\rho] = \sqrt{\sum_{|k|=2}^{K} \rho_k^2}\,, 
\end{equation}
We consider $L^2$ errors of the mean and standard deviation (SD). 

\[
\begin{aligned}
	&\text{Error}_{\text{mean}}(\rho) = \left\| \mathbb{E}[\rho_h] - \mathbb{E}[\rho] \right\|_{L^2}, \quad
	\text{Error}_{\text{SD}}(\rho) = \left\| \mathbb{SD}[\rho_h] - \mathbb{SD}[\rho] \right\|_{L^2}, \\[4pt]
	&\text{Diff}_{\text{mean}}(\rho) = \mathbb{E}[\rho_h] - \mathbb{E}[\rho], \quad
	\text{Diff}_{\text{SD}}(\rho) = \mathbb{SD}[\rho_h] - \mathbb{SD}[\rho],
\end{aligned}
\]
where $\rho_h$ and $\rho$ are respectively numerical solutions of gPC-SG and the reference
solutions, which is obtained by the stochastic collocation (SC) method \cite{Xiu}. Let $\{z^{(j)}\}_{j=1}^{N_z} \subset I_z$ be the set of collocation points, with the corresponding weights $\{w^{(j)}\}$ determined by the quadrature rule and $N_z$ the number of samples. The integral such as expectation of $\rho$ in the random space can be approximated by:
\begin{equation}
\int_{I_z} \rho(t, x, v, z) \pi(z) \, dz \approx \sum_{j=1}^{N_z} \rho(t, x, v, z^{(j)}) w^{(j)}.
\end{equation}
 {In Test I (b) below, we take $N_z=10$ in each dimension of random space, while in other tests  $N_z=16$. } 

\subsection{Test I: Uncertain initial data with finger-like projection property}

In this test, we focus on the early stages of tumor development, particularly the formation of finger-like projections. To solve (\ref{test I and II}), we employ the in vivo model for the nutrient equation \( c(x, z) \). We consider uncertain nutritional function \( c(x, z) \), the initial condition \( f \) characterized by a random amplitude with a finger-like structure, the radius of support of \( \rho \) as well as parametric data \( G_0 \). {We design experiments with one- and two-dimensional random variables and compare the computational cost of SG and SC methods for both 1D and 2D examples, including CPU time and memory tests. Memory usage 
is measured using \texttt{profile(`memory',`on')} in MATLAB. The recorded CPU time and memory values are averaged by multiple simulations to ensure their reliability. }

Set the computational domain \( \Omega = [-2.2, 2.2] \times [-2.2, 2.2] \) with spatial coordinates \( (X, Y) \in \Omega \). Let spatial mesh size \( \Delta x = \Delta y = 0.1 \), with temporal step size \( \Delta t = 0.001 \). Unless otherwise specified, we take \( m = 80 \). 

{\bf{Test I (a).}} The initial condition is given by
\begin{equation}
    \rho (X, Y, z) =
\begin{cases}
0.15 - 0.1z, & r \leq r_1, \\[5pt]
\frac{r - r_1}{5(r_2 - r_1)}+ 0.15 - 0.1z, & r_1 < r \leq r_2, \\[5pt]
0.35 + 0.1z, & r_2 < r \leq r_3, \\[5pt]
\frac{r_3 - r}{r_3 - r_4}(0.35 + 0.1z), & r_3 < r \leq r_4, \\[5pt]
0, & r > r_4.
\end{cases}
\end{equation}

\begin{equation}
\label{growth function for nutrition}
    H(X, Y, c,z) = 
\begin{cases} 
 G_0 c(x,z), & \sqrt{(X )^2 + (Y )^2} \leq r_4, \\
1, & \sqrt{X^2 + Y^2} > r_4.
\end{cases}
\end{equation}
where \( r = \sqrt{X^2 + Y^2} \), and the radii are given by:
\[
r_1 = 0.4(1 + 0.3z), \quad r_2 = 0.5(1 + 0.3z), \quad r_3 = 0.6(1 + 0.3z), \quad r_4 = 0.7(1 + 0.3z).
\]
Let the two nutritional parameters be defined as \( G_0 = 0.5(1 - 0.1z) \). We solve the nutrition matrix \( H \), \( c(x,z) \) from equations (\ref{HHH}) and (\ref{nutrition caculation}). 

\begin{figure}[H]
    \centering
    \includegraphics[width=0.9\textwidth]{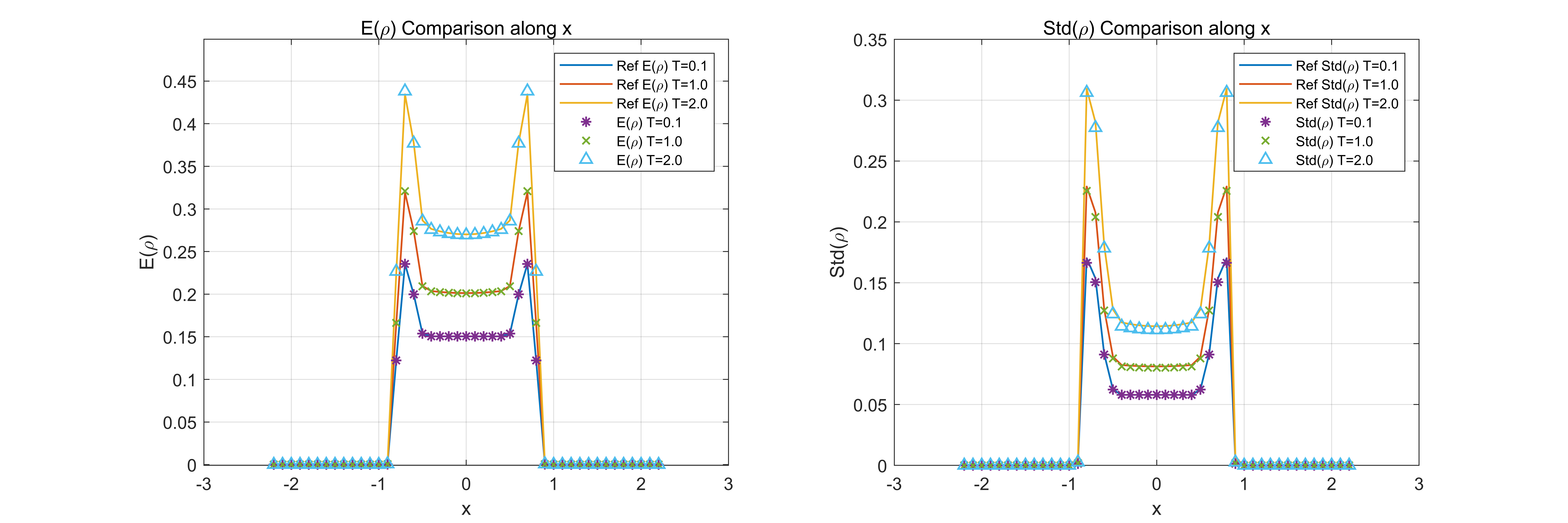}  
\captionsetup{justification=raggedright,singlelinecheck=false}
    \caption{Test I {(a)}: Mean and standard deviation of density $\rho$ (\(x = -0.3\)) at different times. $\Delta x = 0.1$, $\Delta t = 1\times 10^{-3}$
, $m = 80$. Star: gPC-SG with \(K = 4\). Solid line:
reference solutions by SC method using $N_z = 16$.}
    \label{fig:test1_slice}
\end{figure}

We plot the density solution at time $T=1$. 
A satisfactory agreement between the gPC-SG solutions and the reference solutions is observed at different times, as shown in Figure \ref{fig:test1_slice}.  

\begin{figure}[H]
    \centering
    \includegraphics[width=0.6\textwidth]{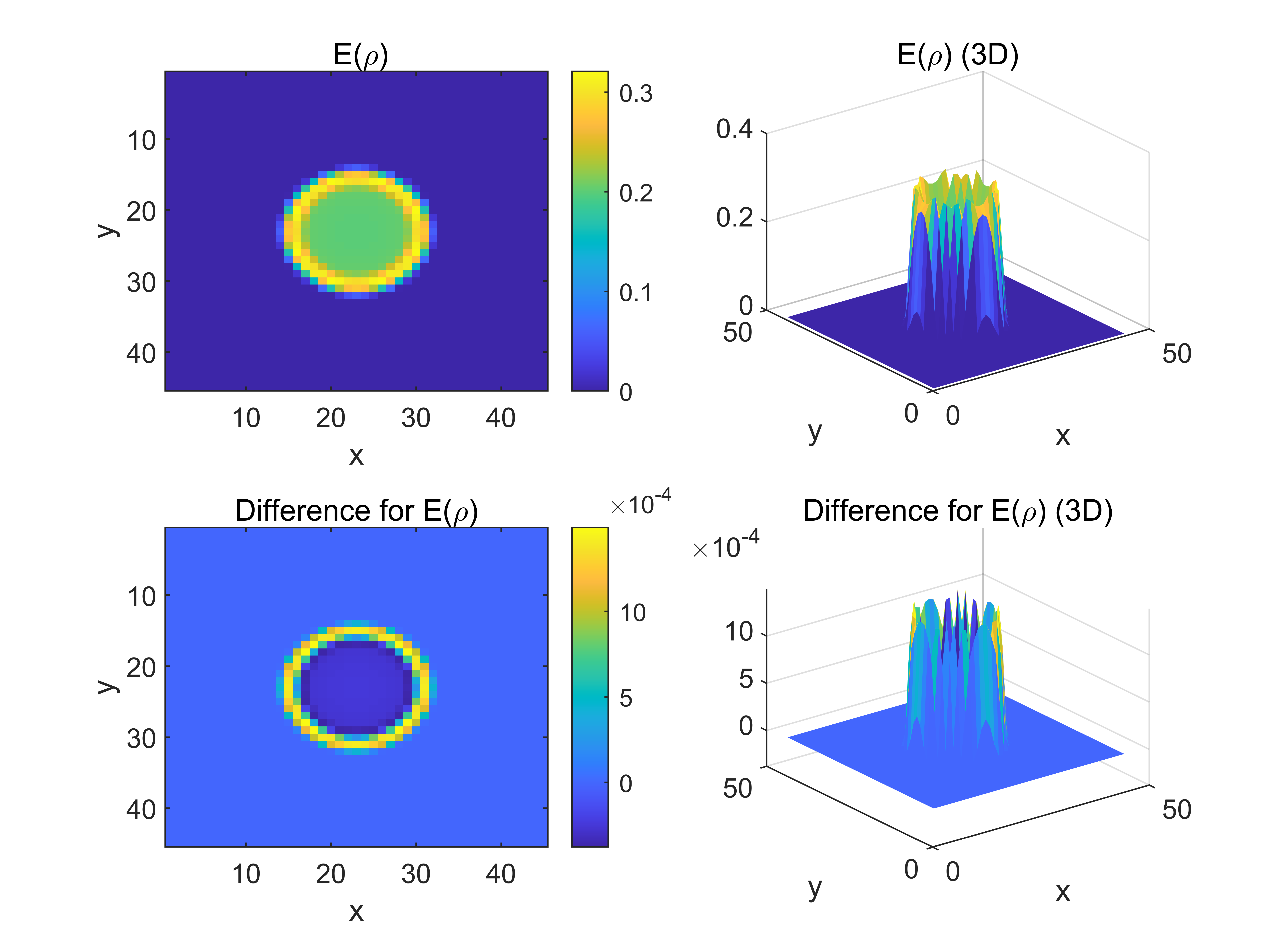}  
\captionsetup{justification=raggedright,singlelinecheck=false}
    \caption{Test I {(a)}: Mean of forth-order gPC-SG solutions at $T = 1$. The second row shows the  differences between the mean (in random space) of gPC-SG and reference solutions at each $x$. }
    \label{fig:test1_mean}
\end{figure}

\begin{figure}[h]
    \centering
    \includegraphics[width=0.6\textwidth]{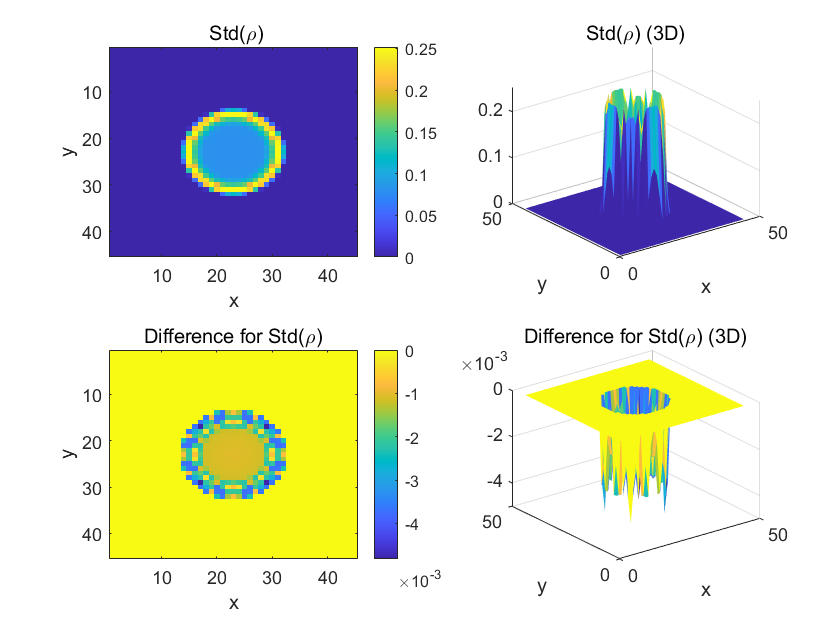}  
\captionsetup{justification=raggedright,singlelinecheck=false}
    \caption{Test I {(a)}: SD of forth-order gPC-SG solutions at $T = 1$. The second row shows the differences between the SD (in random space) of gPC-SG and reference solutions at each $x$. }
    \label{fig:test1_std}
\end{figure}

In Figure \ref{fig:test1_mean}, we present the mean of fourth-order gPC-SG solutions and the differences between the gPC-SG and reference solutions at \( T = 1 \). Figure \ref{fig:test1_std} shows the standard deviation of the gPC-SG solutions and the differences between the two solutions at \( T = 1 \).  

The tumor model exhibits faster growth in regions with higher curvature and is characterized by finger-like projections \cite{cristini2010multiscale,saffman1958penetration}. The porous medium-type model has finger-like properties driven by viscosity variation, and greater curvature of the tip is associated with fuller contact with external nutrients per unit volume, with relative nutrient enrichment and hence faster growth. In Test I (a), the cross-sectional view shown in Figure \ref{fig:test1_slice} reveals that the curvature is larger at the tip of the finger, leading to noticeably rapid growth. Furthermore, we observe that the perturbation has a significant effect in faster-growing regions of high curvature, such as the boundary, where the difference between the gPC-SG solution and the SC solution is more noticeable as can be seen in Figure \ref{fig:test1_mean} and Figure \ref{fig:test1_std}.  The errors in the mean and standard deviation both remain below $\mathcal{O}(10^{-3})$. 

\begin{figure}[H]
    \centering
    \includegraphics[width=0.4\textwidth]{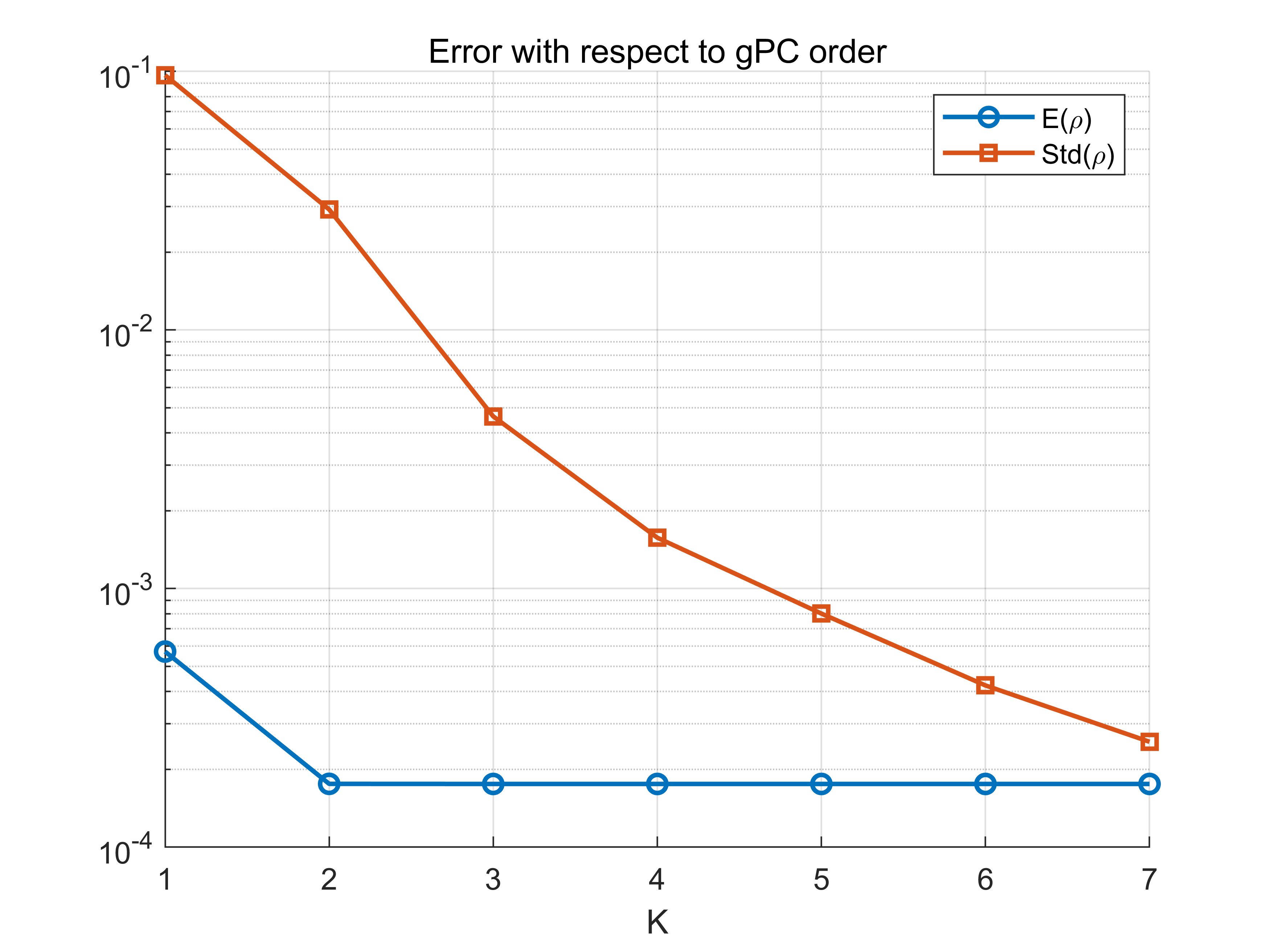}  
    \caption{Test I {(a)}: Errors in space between the gPC-SG and reference solutions for the mean and
SD of $\rho$ with respect to gPC order \(K\) at $T = 0.5, \Delta x = 0.1, \Delta t = 1 \times 10^{-3}$. }
    \label{fig:test1_error}
\end{figure}

Figure \ref{fig:test1_error} illustrates the rapid exponential decay of \( L^2 \) spatial errors for the mean and standard deviation of \( \rho \), for solving the 2D porous medium-type system with varying gPC orders \( K \). We observe that the errors in the mean of \( \rho \) quickly saturate at moderately small gPC orders, particularly when \( K = 2 \). 

{\bf{Test I (b). }}
{
	The initial condition is given by
}
\begin{equation}
	{
		\rho(X, Y, z_1, z_2) =
		\begin{cases}
			0.15 - 0.1 z_1, & r < r_1, \\[5pt]
			\dfrac{r - r_1}{10(r_2 - r_1)}  (2 + z_2+ z_1) + (0.15 - 0.1 z_1), & r_1 \leq r < r_2, \\[5pt]
			0.35 + 0.1 z_2, & r_2 \leq r < r_3, \\[5pt]
			\left( \dfrac{r_4 - r}{r_4 - r_3} \right)(0.35 + 0.1 z_2), & r_3 \leq r < r_4, \\[5pt]
			0, & r \geq r_4,
		\end{cases}
	}
\end{equation}
{
	where \( z_1 \) and \( z_2 \) are independent random variables uniformly distributed on the interval \([-1, 1]\). Define \( r = \sqrt{X^2 + Y^2} \), 
    we assume the radii contain random parameters which are given by
}
\begin{equation}
	{
		r_1 = 0.4 (1 + 0.2 z_1), \quad
		r_2 = 0.5 (1 + 0.2 z_1), \quad
		r_3 = 0.6 (1 + 0.2 z_2), \quad
		r_4 = 0.7 (1 + 0.2 z_2),
	}
\end{equation}

{
	The growth function for nutrient consumption is set by
}
\begin{equation}
	{
		\label{growth_function}
		H(X, Y, c, z_2) =
		\begin{cases}
			G_0 \, c(X, Y, z_2), & \sqrt{X^2 + Y^2} \leq r_4, \\[5pt]
			1, & \sqrt{X^2 + Y^2} > r_4,
		\end{cases}
	}
\end{equation}
{
	where \( G_0 = 0.5(1 - 0.1z_1- 0.1z_2) \), and \( c(X, Y, z_2) \) is the nutrient concentration. 
    }

\begin{figure}[H]
    \centering
    \includegraphics[width=0.9\textwidth]{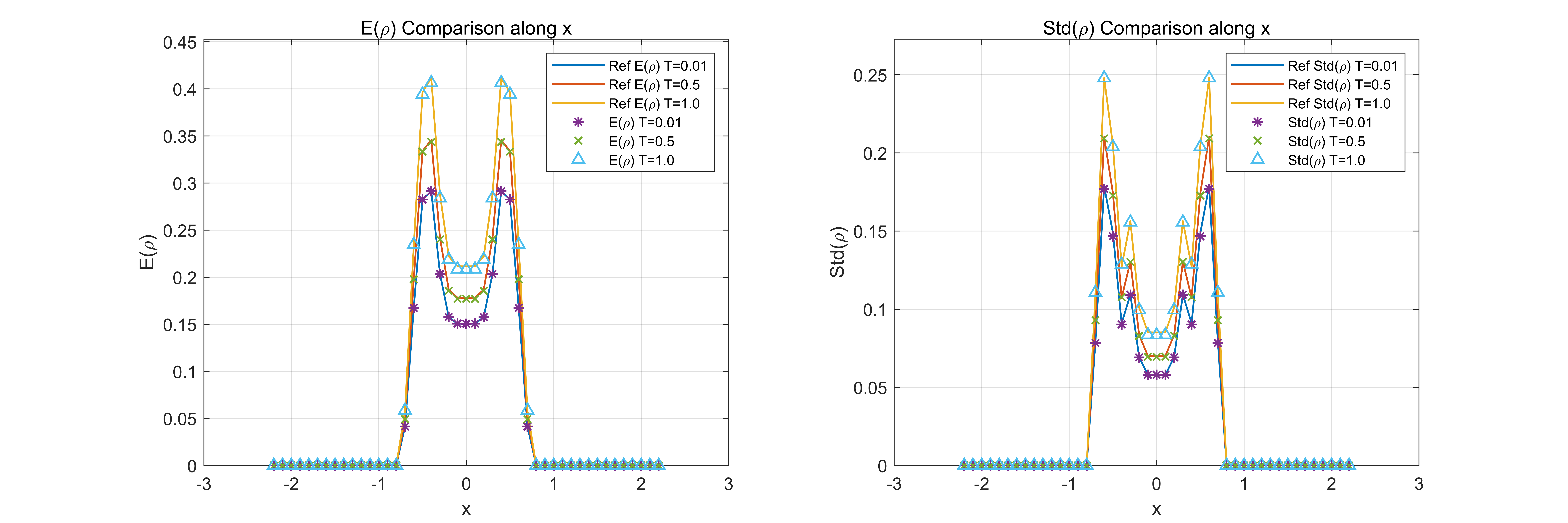}  
\captionsetup{justification=raggedright,singlelinecheck=false}
    \caption{{Test I (b): Mean and standard deviation of density $\rho$ (\(x = -0.3\)) at different times. $\Delta x = 0.1$, $\Delta t = 1\times 10^{-3}$
, $m = 80$. Star: gPC-SG with \(K = 4\). Solid line:
reference solutions by SC method using $N_z = 10$ in each dimension.}}
\label{fig:testI_b_slice}
\end{figure}

{In Figure~\ref{fig:testI_b_slice}, a satisfactory agreement between the gPC-SG solutions and the reference solutions is observed at different times. Compared to the one-dimensional in random variable case, SG method is still able to capture detailed characteristics of the solutions, as shown in the plots of standard deviation. }


\begin{figure}[h]
    \centering
    \includegraphics[width=0.6\textwidth]{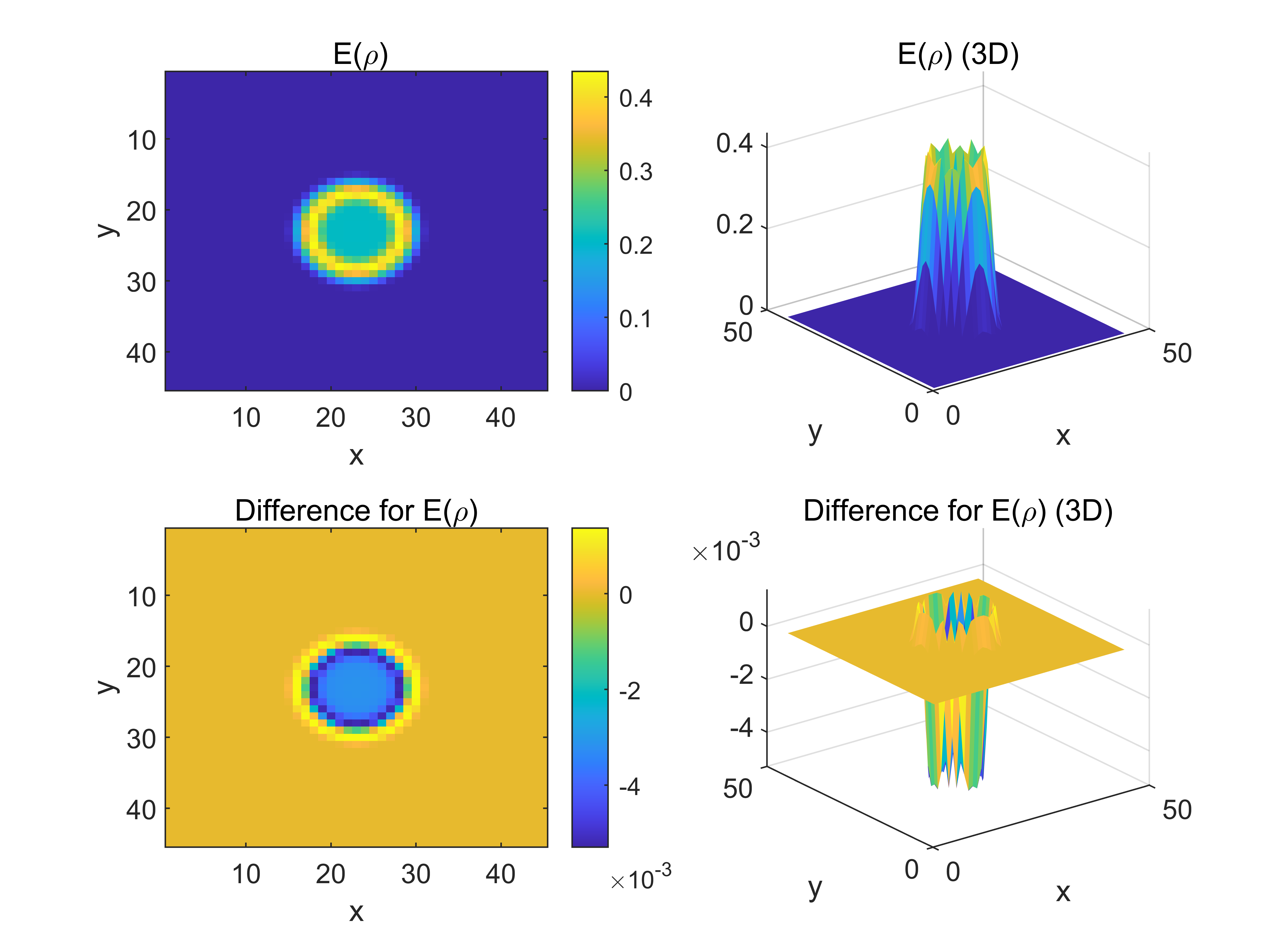}  
\captionsetup{justification=raggedright,singlelinecheck=false}
    \caption{{Test I (b): Mean of forth-order gPC-SG solutions at $T = 1$. The second row shows the differences between the mean (in random space) of gPC-SG and reference solutions at each $x$.}}
    \label{fig:testI_b_mean}
\end{figure}

\begin{figure}[H]
    \centering
    \includegraphics[width=0.6\textwidth]{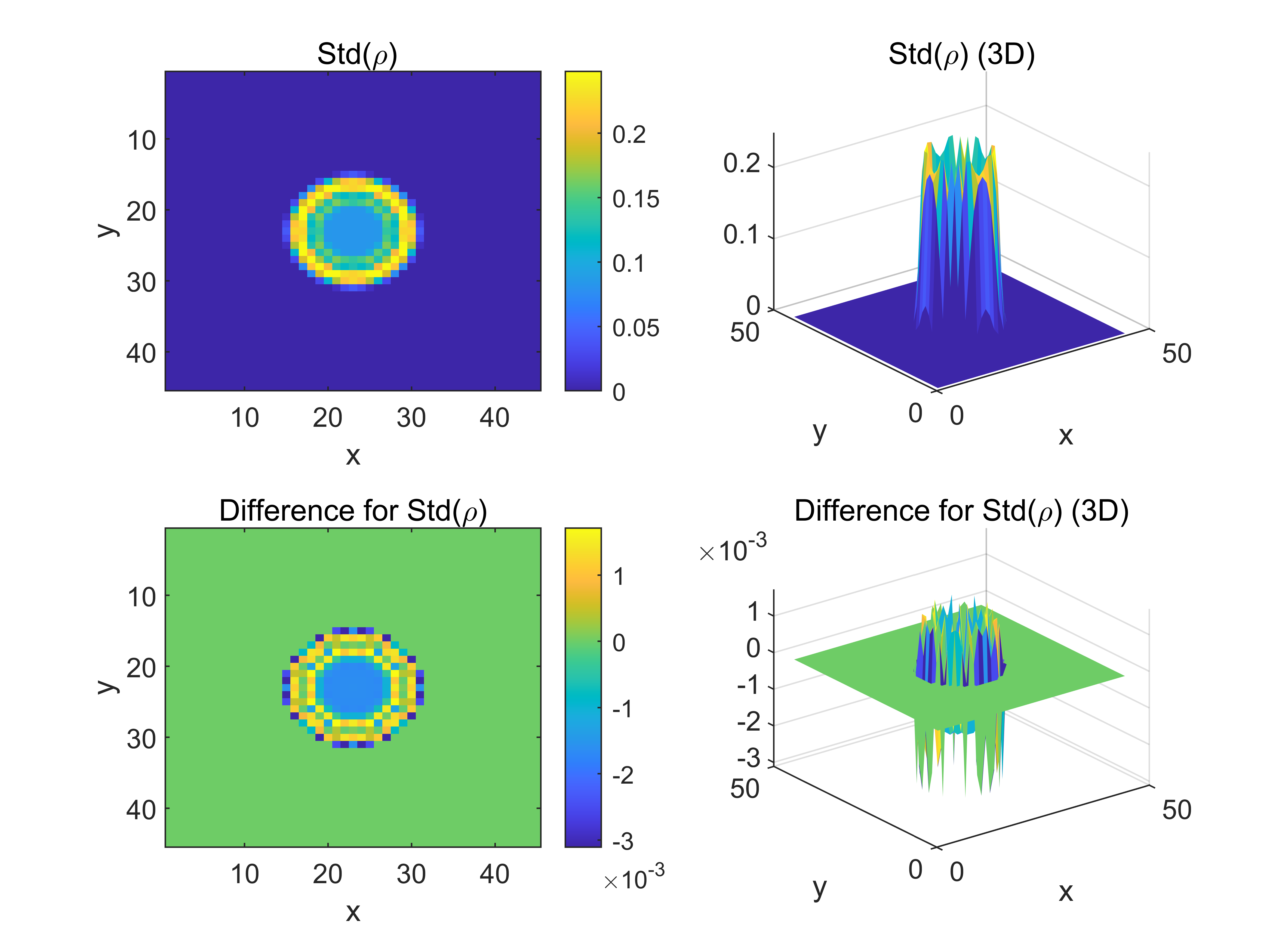}  
\captionsetup{justification=raggedright,singlelinecheck=false}
    \caption{{Test I (b): SD of forth-order gPC-SG solutions at $T = 1$. The second row shows the differences between the SD (in random space) of
    gPC-SG and reference solutions at each $x$.}}
    \label{fig:testI_b_std}
\end{figure}

{In Figure~\ref{fig:testI_b_mean}, we present the mean of the fourth-order gPC-SG solutions, along with the differences between the gPC-SG and reference solutions at} {$T = 1$}{. Figure~\ref{fig:testI_b_std} shows the similar plots for the standard deviation of gPC-SG solutions. }

\begin{figure}[H]
    \centering
    \includegraphics[width=0.4\textwidth]{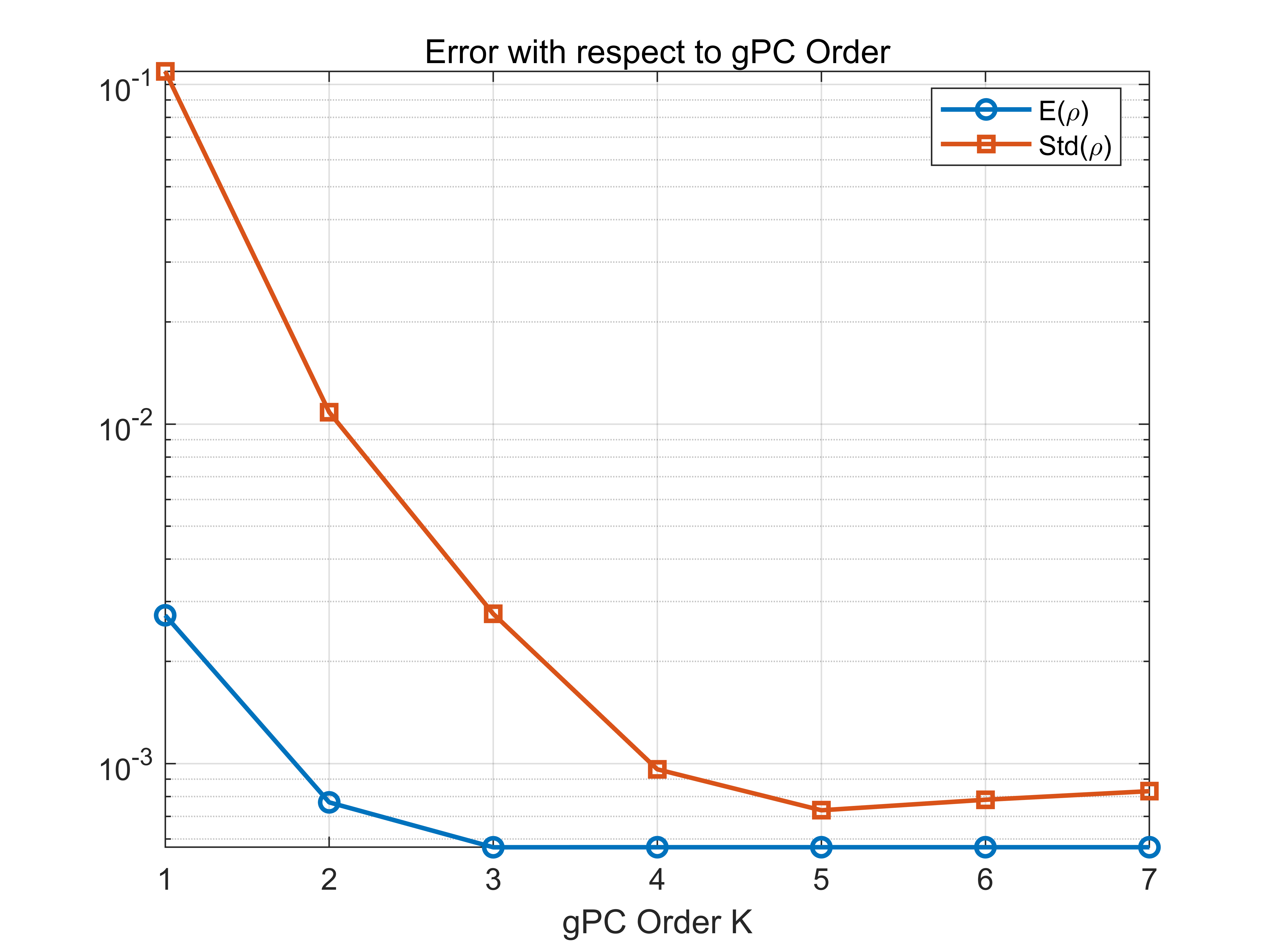}  
    \caption{{{Test I (b): Errors in space between the gPC-SG and reference solutions for the mean and
SD of $\rho$ with respect to gPC order $K$ at
$T = 0.5$, $\Delta x = 0.1$, 
$\Delta t = 1 \times 10^{-3}$. }}}
    \label{fig:test1(b)_error}
\end{figure}

{Figure~\ref{fig:test1(b)_error} illustrates a rapid exponential decay of } {$L^2$} {errors in both mean and standard deviation of the density solution, as the gPC order} {$K$} {increases. Notably, the errors of mean saturate as small as} {$K = 3$}{, indicating a fast spectral convergence. This behavior is consistent with that observed in Test I (a).}

{In Table~\ref{cpu time}, we demonstrate that SG method significantly outperforms SC method in both 1D and 2D random variable problems. In Test 1 (a), SG method is over 20 times faster than SC approach; in Test 1 (b), SG method runs about 3 times faster. More importantly, SG maintains this advantage in computational speed as the dimensionality of random parameter increases. This consistent performance in CPU time highlights the efficiency and robustness of SG method.}

\begin{figure}[H]
    \centering
    \includegraphics[width=0.4\textwidth]{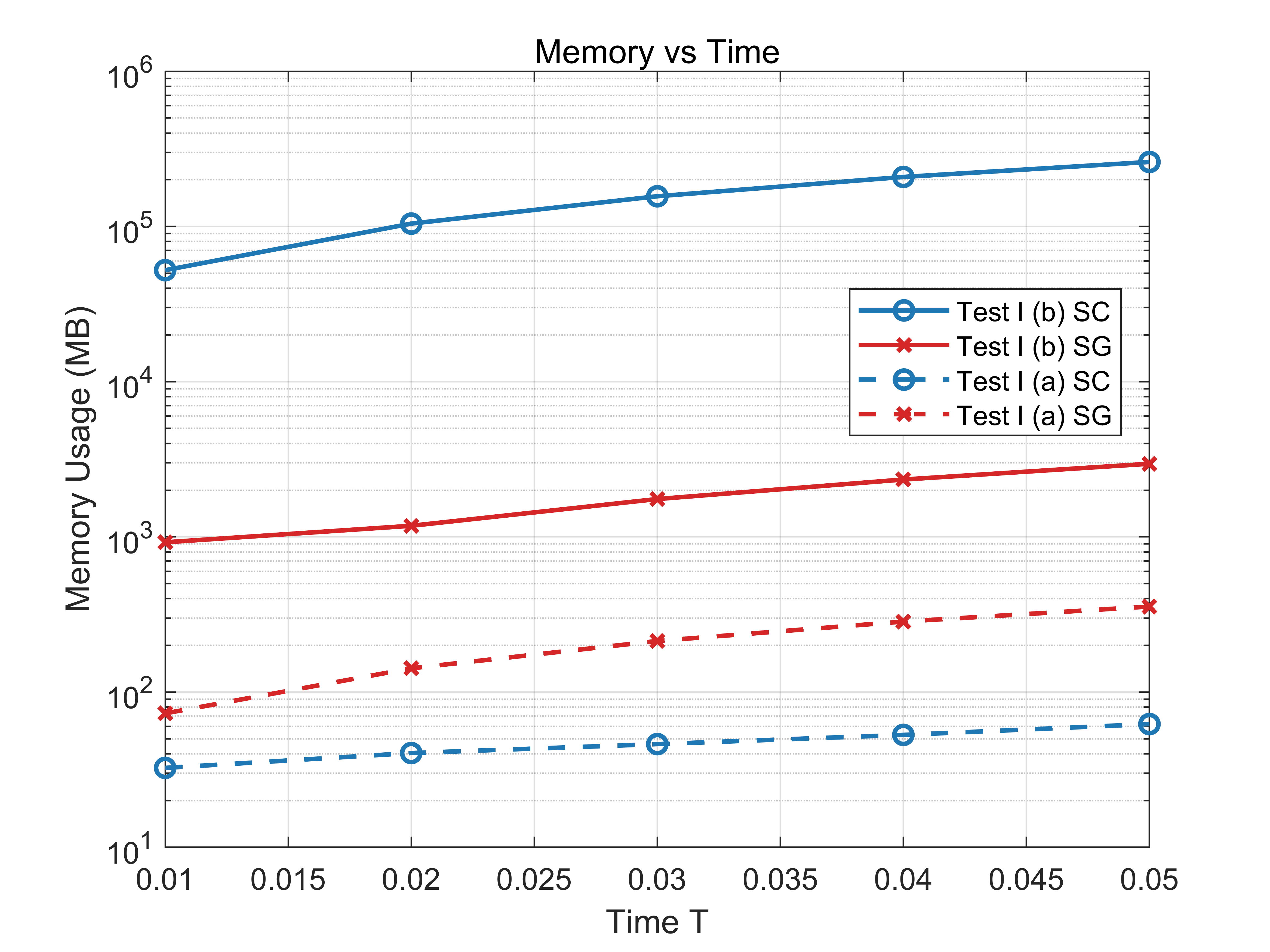}  
    \caption{{Test I: memory of SG and SC at different times for Test I. }}
    \label{fig:memory_vs_time}
\end{figure}

{{Besides the CPU time, a memory usage comparison is shown in Figure~\ref{fig:memory_vs_time}, which further validates the efficiency of SG method. For the one-dimensional problem in Test I (a), both SC and SG methods exhibit relatively low memory consumption--approximately $O(10^2)$MB for SC and $O(10^3)$MB for SG, showing only a slight difference. However, as the dimensionality of random variables increases in Test I (b), the gap becomes dramatic and significant: SC consistently consumes over 50 times more memory than SG for different output times, and nearly 90 times more at the final time $T=0.05$. This phenomenon further highlights the excellent scalability and memory advantages of the SG over SC method for high-dimensional problems. }}
    
\begin{table}[htbp]
	\centering
	\begin{tabular}{ccccc}
		\hline
		{Test} & {Method} & {$T = 0.01$} & {$T = 0.03$} & {$T = 0.05$} \\
		\hline
		{Test I (a)} 
		& {SG}  & {0.5696} & {1.3974} & {2.3244} \\
		& {SC}  & {10.9827} & {30.1244} & {48.1486} \\
		{Test I (b)} 
		& {SG}  & {6.3211} & {17.9350} & {29.6033} \\
		& {SC}  & {18.4569} & {54.1962} & {89.4438} \\
		\hline
	\end{tabular}
    	\caption{{CPU time (units in seconds) of SG and SC at different times for Test I.}}
	\label{cpu time}
\end{table}

\subsection{Test II: Uncertain initial data with general planar features}

In this test, we examine early-stage tumor growth with general planar structures and adopt the in vivo nutrition model by (\ref{test I and II}), (\ref{nutrition caculation}). The sources of randomness involve the nutritional function \( c(x,z) \), the initial condition \( f \) characterized by a random amplitude with general planar features, the radius of support of \( \rho \) and the parametric data \( G_0 \). The initial condition is given by
\begin{equation}
    \rho(X, Y, z) =
\begin{cases} 
0.2(1 + 0.5z), &  \sqrt{X^2 + Y^2} \leq r_1, \\[5pt]
\frac{r_2 - \sqrt{X^2 + Y^2}}{5(r_2 - r_1)}(1 + 0.5z), &  r_1 < \sqrt{X^2 + Y^2} \leq r_2, \\[5pt]
0, &  \sqrt{X^2 + Y^2} > r_2,
\end{cases}
\end{equation}

\begin{equation}
    H(X, Y, z,c) = 
\begin{cases} 
 G_0c(x,z), & \sqrt{(X )^2 + (Y )^2} \leq r_2, \\
1, & \sqrt{X^2 + Y^2} > r_2,
\end{cases}
\end{equation}
where \[
r_1 = 0.4(1 + 0.4z), \quad r_2 = 0.5(1 + 0.5z).
\]

\begin{figure}[h]
    \centering
    \includegraphics[width=0.9\textwidth]{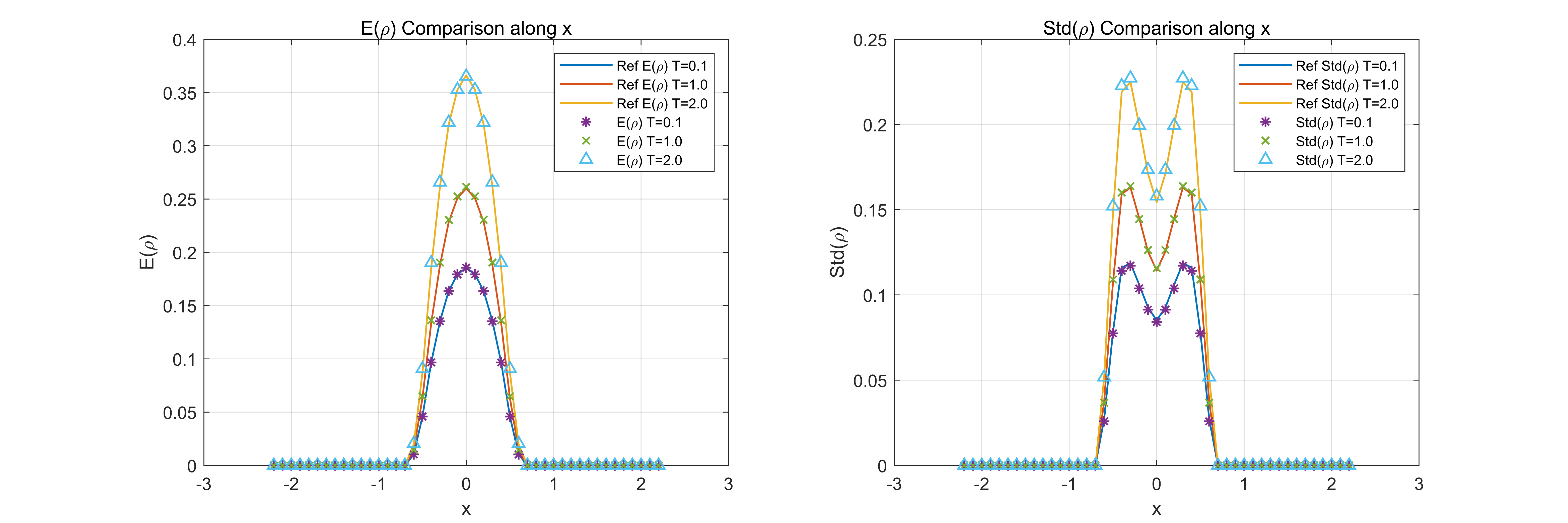}   \captionsetup{justification=raggedright,singlelinecheck=false}
    \caption{Test II: Mean and standard deviation of density $\rho$ (\(x = -0.3\)) at different time \(T\). $\Delta x = 0.1$, $\Delta t = 1\times 10^{-3}$
, \(m = 80\). Star: gPC-SG with \(K = 4\). Solid line:
the reference solutions by the collocation method using \(N_z = 16\).}
    \label{fig:testII-slice}
\end{figure}

The growth rate function \( H(X,Y,c,z) \) remains the same as in Test I, see (\ref{growth function for nutrition}). We employ the SC method to obtain reference solutions. A comparison between the gPC-SG and reference solutions at different times is presented in Figure \ref{fig:testII-slice}, which demonstrates a satisfactory agreement between the two solutions.  

\begin{figure}[tpb]
    \centering
    \includegraphics[width=0.6\textwidth]{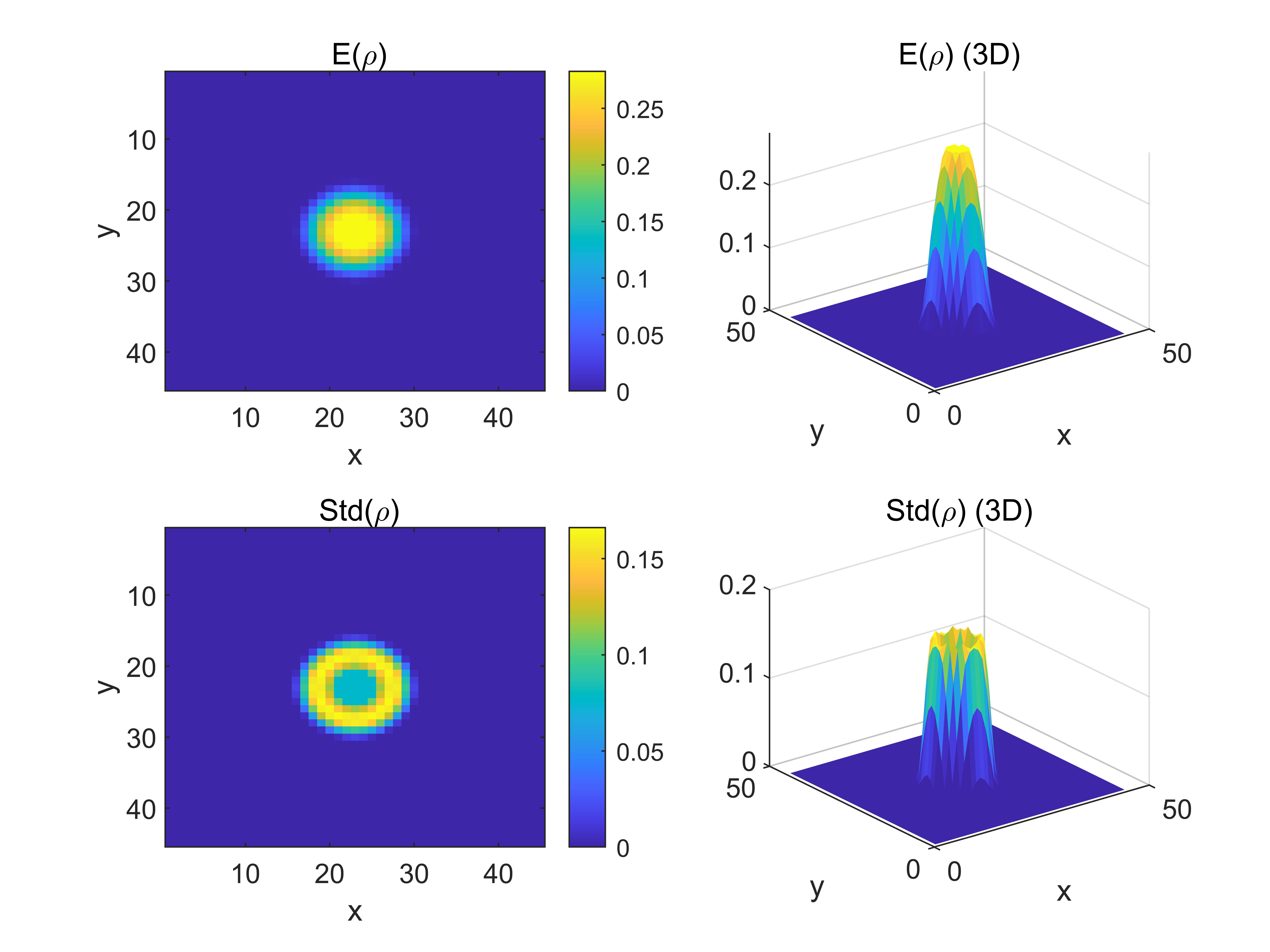}  \captionsetup{justification=raggedright,singlelinecheck=false}
    \caption{Test II. First row:  Mean of the forth-order gPC-SG solutions at \(T = 1\). Second row: Standard deviation of the forth-order gPC-SG solutions at \(T = 1\).}
    \label{fig:testII. mean}
\end{figure}

In Figure \ref{fig:testII. mean}, we present the mean and standard deviation of fourth-order gPC-SG solution at time \( T = 1 \). In Test II, regions with higher central curvature grow faster. This phenomenon is clearly reflected in the section diagram seen in Figure \ref{fig:testII-slice}. This behavior is consistent with the basic principle behind the finger-like structure observed in Test I, particularly in areas of higher curvature, the contact area per unit volume is larger, nutrient-enriched, and the growth rate is faster.

\subsection{Test III: Uncertain parameters in the proliferating, quiescent and dead cells model}

In this test, we examine a tumor growth model that incorporates proliferating, quiescent and dead cells. Consider in vivo nutrition model, 
sources of uncertainty come from the nutritional function \( c(x,z) \), initial conditions of \( \rho_P \), \( \rho_Q \), and \( \rho_D \) with random amplitude, radius of support of the tumor \( \rho \) as well as parametric data \( a \), \( b \), \( d \), and \( \mu \). Let \( \rho_P(x,t) \), \( \rho_Q(x,t) \), and \( \rho_D(x,t) \) represent the cell densities of proliferating, quiescent, and dead cells, respectively. The total cell density is given by
\[
\rho(x,t) = \rho_P(x,t) + \rho_Q(x,t) + \rho_D(x,t).
\]
We consider the initial data as follows
\begin{equation}
    \rho(X, Y, z) =
	\begin{cases}
		0.3 (1+ 0.5z), & r \leq r_1, \\[5pt]
		\frac{r - r_1}{r_2 - r_1}(0.3(1+ 0.5z)) + (0.3(1+ 0.5z)), & r_1 < r \leq r_2,
	\end{cases}
\end{equation}

\begin{equation}
    \rho_D(X, Y, z) =
	\begin{cases}
		0.15 + 0.1z, & r \leq r_3, \\[5pt]
		\frac{r - r_3}{r_4 - r_3}(0.15 + 0.1z) + (0.15 + 0.1z), & r_3 < r \leq r_4,
	\end{cases}
\end{equation}
and
\begin{equation}
    \rho_P(X, Y, z) = 0.55 \left( \rho(X, Y, z) - \rho_D(X, Y, z) \right),
\end{equation}
\begin{equation}
    \rho_Q(X, Y, z) = 0.45 \left( \rho(X, Y, z) - \rho_D(X, Y, z) \right),
\end{equation}
where  \( r = \sqrt{X^2 + Y^2 - 0.8 \sin(4 \, \text{atan2}(Y, X))} \). The radii are given by:	
	\[
	r_1 = 0.4(1 + 0.5z), \quad r_2 = 0.5(1 + 0.5z)\quad
	  r_3 = 0.2(1 + 0.2z), \quad r_4 = 0.3(1 + 0.4z). \]
\[a = 0.2(1+0.25z);\quad
b = 0.2(1-0.25z);\quad
d = 0.25(1-0.6z);\quad
\mu = 0.\]
       
\begin{figure}[H]
    \centering
    \includegraphics[width=0.9\textwidth]{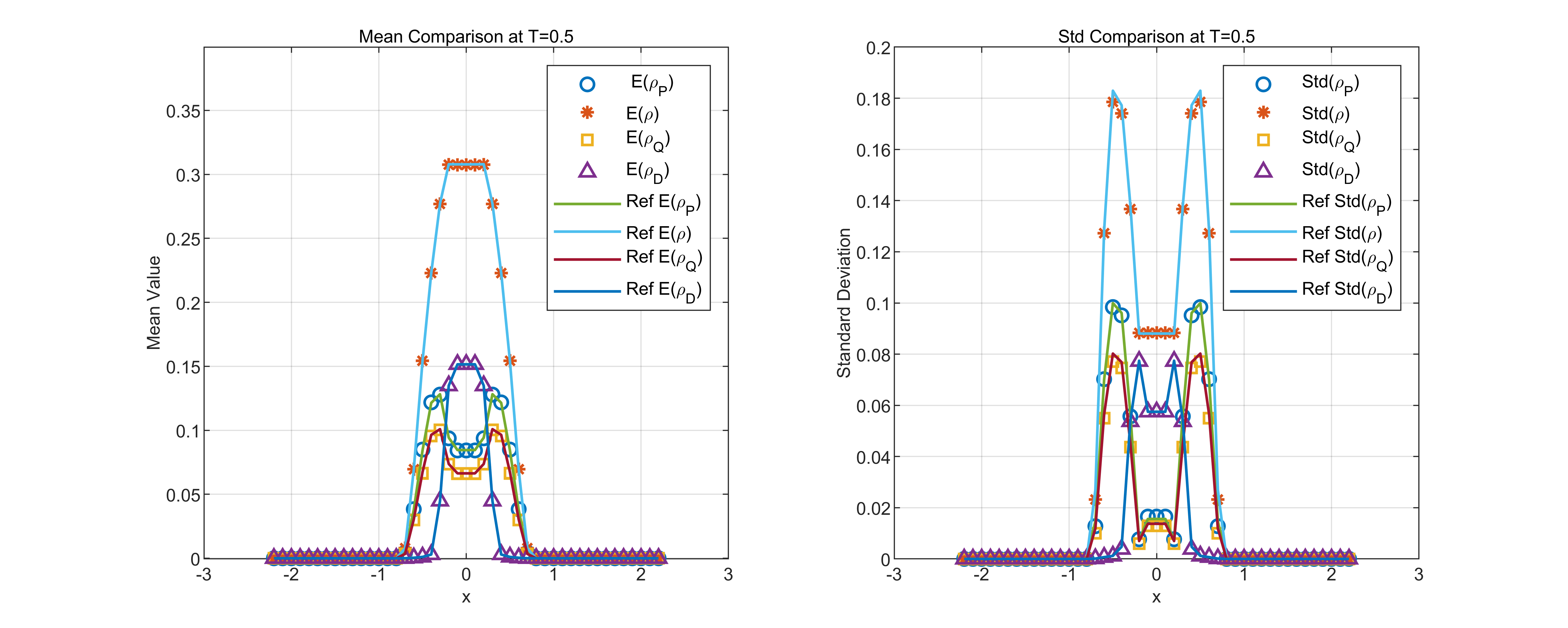}  
    \caption{Test III: Mean and standard deviation of density $\rho$, $\rho_P$, $\rho_Q$ and $\rho_D$ at \(T = 0.5\), \(x = 0\). $\Delta x = 0.1$, $\Delta t = 1\times 10^{-3}$
, \(m = 80\). Star: gPC-SG solution with \(K = 4\). Solid line: reference solutions by SC method using \(N_z = 16\).}
    \label{fig:test3_slice}
\end{figure}

We obtain similar results: a satisfactory agreement between the gPC-SG solutions and the reference solutions is observed for \( \rho \), \( \rho_P \), \( \rho_Q \), and \( \rho_D \), as shown in Figure \ref{fig:test3_slice}.  

\begin{figure}[H]
    \centering
    \includegraphics[width=0.5\textwidth]{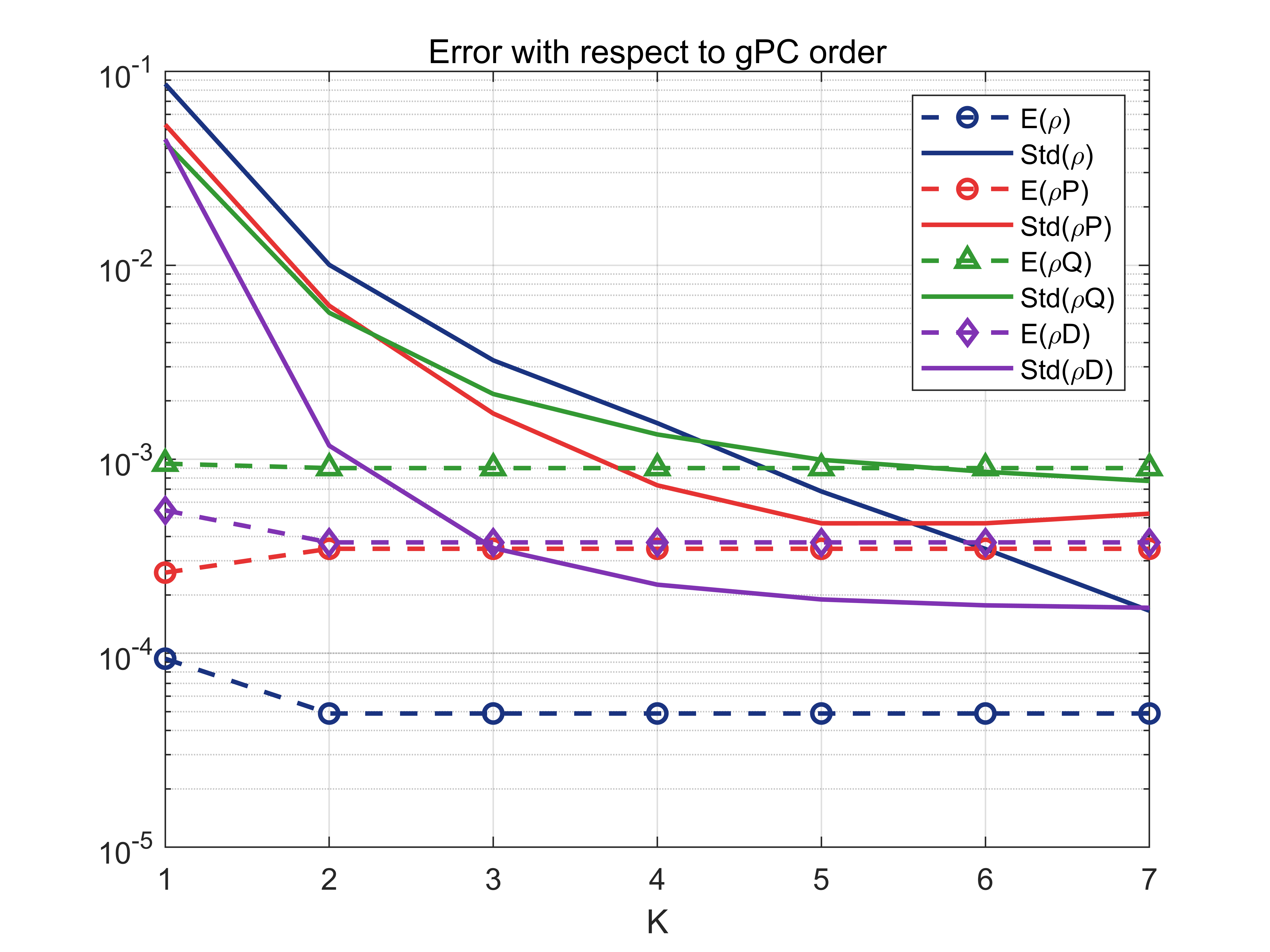}  
    \caption{Test III: $L^2$ errors between the gPC-SG and reference solutions for the mean and
standard deviation of $\rho$, $\rho_P$, $\rho_Q$ and $\rho_D$ with respect to different $K$ at $T = 0.1$. Set $\Delta x = 0.1, \Delta t = 5 \times 10^{-3}$. }
    \label{fig:test3_error}
\end{figure}

Figure \ref{fig:test3_error} shows the rapid exponential decay of \( L^2 \) errors for the standard deviation of \( \rho_P \), \( \rho_Q \), and \( \rho_D \) with respect to different gPC orders \( K \). We observe that the errors in the mean values of \( \rho \) and \( \rho_D \) quickly saturate at $K=2$. 

\begin{figure}[htpb]
    \centering
    \includegraphics[width=0.6\textwidth]{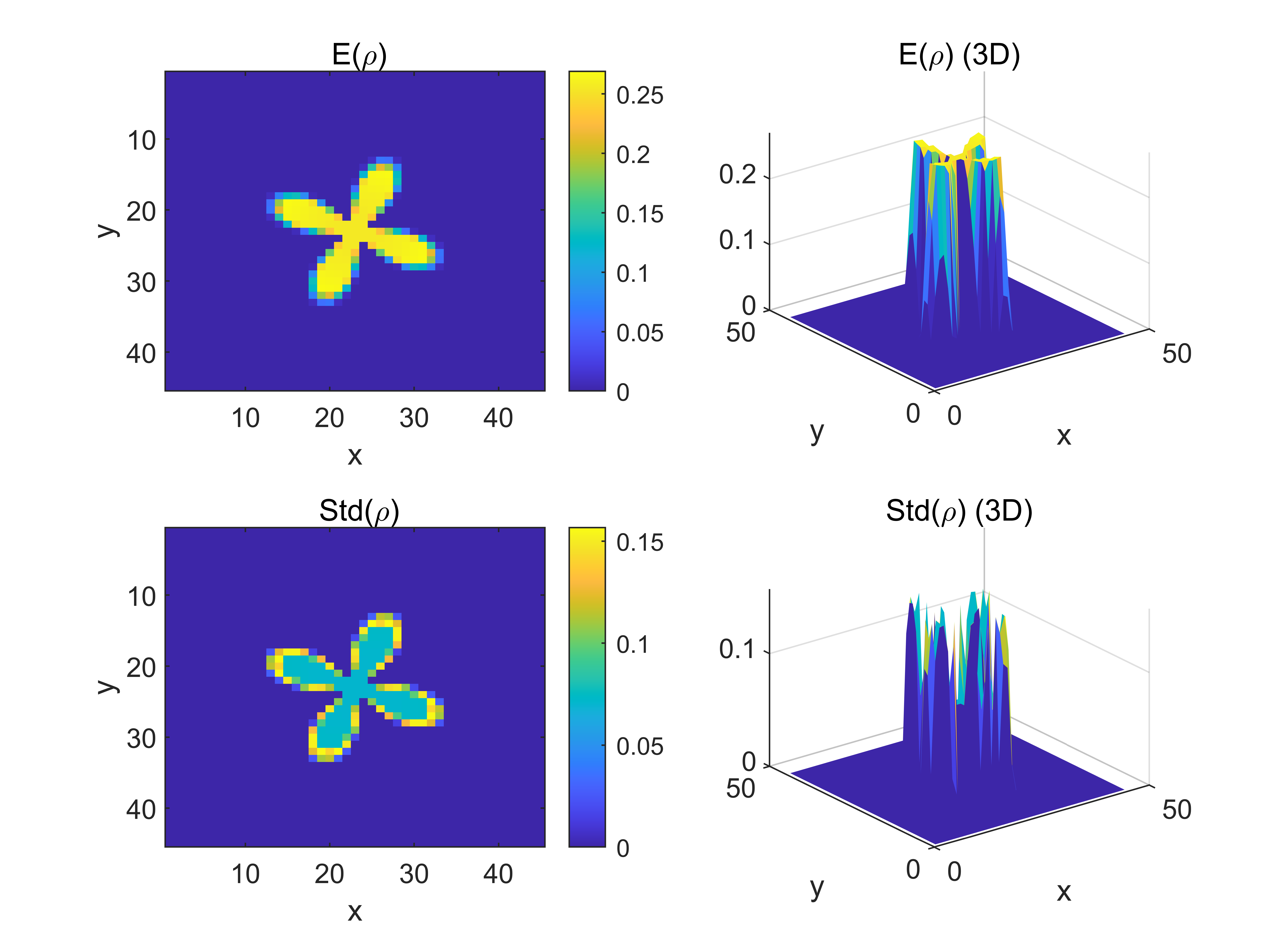} \captionsetup{justification=raggedright,singlelinecheck=false}
    \caption{Test III: Mean of the forth-order gPC-SG solutions of $\rho$ at $T = 0.5$. The second row: standard deviation of the forth-order gPC-SG solutions of $\rho$ at $T = 0.5$.}
    \label{fig:test3_mean}
\end{figure}

\begin{figure}[htpb]
    \centering
    \includegraphics[width=0.6\textwidth]{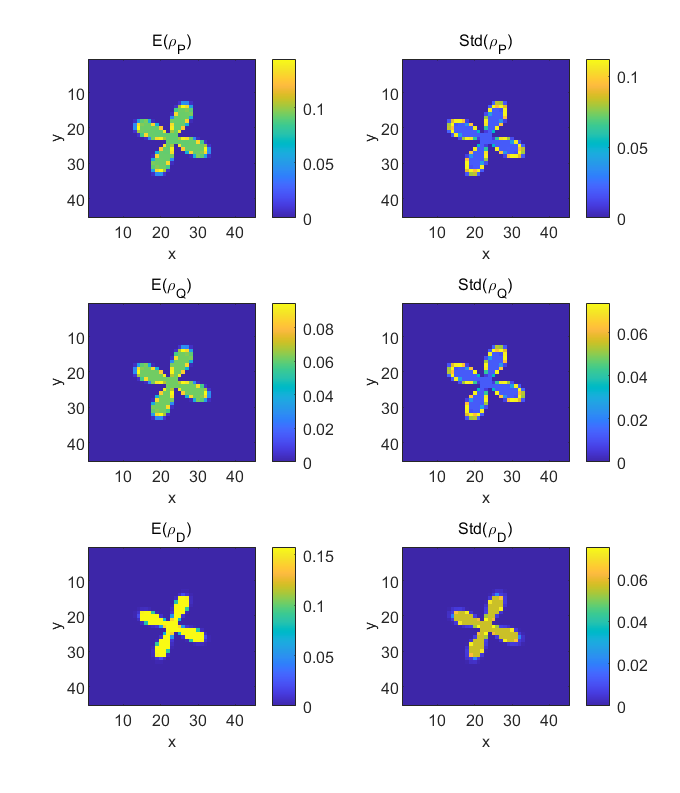}  
\captionsetup{justification=raggedright,singlelinecheck=false}
    \caption{Test III: Mean of the forth-order gPC-SG solutions of $\rho_ P$, $\rho _Q$ and $\rho _D$ at T = 0.5 (left). The second row: standard deviation of the forth-order gPC-SG solutions of $\rho _P$, $\rho _Q$ and $\rho _D$ at T = 0.5 ({right}).}
    \label{fig:test3_3mean}
\end{figure}

In Figure \ref{fig:test3_mean}, we show the mean and standard deviation of fourth-order gPC-SG solutions, and \( L^2 \) errors between the gPC-SG and reference solutions for \( \rho \) at time \( T = 0.5 \). Figure \ref{fig:test3_3mean} shows the mean and standard deviation of gPC-SG solutions for \( \rho_P \), \( \rho_Q \), and \( \rho_D \) at time \( T = 0.5 \). For \( \rho \), \( \rho_P \), \( \rho_Q \), and \( \rho_D \), with more complex boundaries, the characteristics of the variation in mean and standard deviation can still be captured by the gPC-SG method. 
 
As shown in Figure \ref{fig:test3_error}, for the mean and standard deviation, the corresponding error orders of \( \rho_P \), \( \rho_Q \), and \( \rho_D \) are around or lower than \( O(10^{-3}) \).  

\section{Conclusion}

In this paper, we propose a stochastic asymptotic-preserving (s-AP) scheme in the framework of gPC-SG method for solving porous medium-type tumor growth models with uncertainties. The s-AP scheme is designed based on \cite{liu2018accurate} for deterministic problems, which can efficiently capture tumor interface evolutions and pattern formations while remaining stable as the index parameter $m \to \infty$. On the other hand, we analyze the regularity of the solution in the random space and establish the s-AP property on the continuous level. We demonstrate that the SG system converges to the SG system of Hele-Shaw dynamics as \(m \to \infty\). 

A series of numerical experiments, {in both 1D and 2D random variable cases}, are conducted to investigate the effect of uncertainties on tumor progression, and to show the finger-like projections and general planar features in proliferating, quiescent or dead cells. The accuracy and efficiency of our proposed numerical method is validated. Our solutions can capture the morphological changing features of early stage tumors, including the finger-like projections, the case with complex boundaries or multiple growth functions. These findings hopefully can provide some new perspectives on uncertainty quantification modeling and tumor growth dynamics. 
In the future, we will study more complex models that incorporate biological behaviors such as nutrient transport and mechanical interactions. In addition, we will study higher-dimensional UQ problems and develop new and more  advanced numerical methods.

\appendixpage

\textbf{A. A time-splitting method based
	on prediction-correction}

For readers convenience, we briefly review the time-splitting method introduced in \cite{liu2018accurate}. Consider 
\begin{equation}
	\begin{cases}
		\partial_t u = m\nabla (\rho^{m-2}) \nabla \cdot (\rho u) - \rho G(c), \\[4pt]
		\partial_t \rho + \nabla \cdot (\rho u) = \rho G(c), \\[4pt]
		u(x, 0) = -\frac{m}{m - 1} \nabla \rho^{m-1}.
	\end{cases}
	\label{pre-cor}
\end{equation}
We introduce the discrepancy term
\begin{equation}
	W = u + \frac{m}{m - 1} \nabla \rho^{m-1}. 
\end{equation}

To illustrate the propagation of the discrepancy when numerical error is present, one adds a small perturbation to (\ref{pre-cor})
\begin{equation}
	\begin{cases}
		\partial_t u = m\nabla (\rho^{m-2}) \nabla \cdot (\rho u) - \rho G(c) + \delta_1, \\[4pt]
		\partial_t \rho + \nabla \cdot (\rho u) = \rho G(c) + \delta_2, \\[4pt]
		u(x, 0) = -\frac{m}{m - 1} \nabla \rho^{m-1}.
	\end{cases}
\end{equation}
Here \(\delta_1, \delta_2\) are small perturbation functions. One considers
\begin{equation}
	\partial_t W (x,t) = \delta_1 + m\nabla (\rho^{m-2} \delta_2).
\end{equation}

In particular, we study the following equation for the discrepancy: 
\begin{equation}
	\partial_t W (x,t) = - \frac{1}{\epsilon^2} W (x,t),
\end{equation}
where \(0 < \epsilon \leq 1\) is the relaxation constant. This discrepancy equation leads to the following relaxation system
\begin{equation}
	\begin{cases}
		\partial_t u = m\nabla (\rho^{m-2}) \nabla \cdot (\rho u) - \rho G(c) - \frac{1}{\epsilon^2} \left( u + \frac{m}{m - 1} \nabla \rho^{m-1} \right), \\[4pt]
		\partial_t \rho + \nabla \cdot (\rho u) = \rho G(c), \\[4pt]
		u(x, 0) = -\frac{m}{m - 1} \nabla \rho^{m-1}.
	\end{cases}
	\label{pre22}
\end{equation}
The following time-splitting method is adopted: 
\begin{equation}
	\begin{cases}
		\partial_t \rho + \nabla \cdot (\rho u) = \rho G(c), \\[4pt]
		\partial_t u = m\nabla (\rho^{m-2}) \nabla \cdot (\rho u) - \rho G(c),
	\end{cases}
\end{equation}
\begin{equation}
	\begin{cases}
		\partial_t \rho = 0, \\[4pt]
		\partial_t u = - \frac{1}{\epsilon^2} \left( u + \frac{m}{m - 1} \nabla \rho^{m-1} \right).
	\end{cases}
\end{equation}

\textbf{B. Proof of Theorem 2}
\\
\textit{proof:}	 Similar to the proof of Theorem 2 in \cite{benilan1996singular}, for \( m>1 \) the operator is defined by
\begin{equation}\label{eq:scalarOp}
	A_m \rho = \operatorname{sign}(\rho)\,|\rho|^{m-1}\Delta \rho,\quad 
	D(A_m) = \{ \rho\in L^m(\Omega) : \; |\rho|^{m-1}\rho \in D(L) \},
\end{equation}
where the Dirichlet--Laplace operator \( L \) is given by
\begin{equation}\label{eq:Laplace}
	L\rho = \Delta \rho,\quad D(L) = \{ \rho\in W^{1,1}_0(\Omega) : \; \Delta \rho\in L^1(\Omega) \}.
\end{equation}

For the vectorized system, each scalar operator becomes a block operator acting on the vector \( \boldsymbol{\rho}(t,x) \in \mathcal{X} \). We denote the vectorized operator by \( \mathbf{A}_m \), then \( \mathbf{A}_m \) is expressed by
\begin{equation}\label{eq:vecOp}
	\mathbf{A}_m \,\boldsymbol{\rho} \;=\; 
	\begin{pmatrix}
		\bigl( A_m \rho_0 \bigr)\\[1mm]
		\bigl( A_m\rho_1 \bigr)\\[1mm]
		\vdots\\[1mm]
		\bigl( A_m \rho_K \bigr)
	\end{pmatrix}
	\;=\; 
	\begin{pmatrix}
		\bigl( \operatorname{sign}(\rho_0) \, |\rho_0|^{m-1}\Delta \rho_0 \bigr)\\[1mm]
		\bigl( \operatorname{sign}(\rho_1) \, |\rho_1|^{m-1}\Delta \rho_1 \bigr)\\[1mm]
		\vdots\\[1mm]
		\bigl( \operatorname{sign}(\rho_K) \, |\rho_K|^{m-1}\Delta \rho_K \bigr)
	\end{pmatrix}.
\end{equation}

Similarly, for \( m=\infty \), the multivalued operator \( \mathbf{A}_\infty \) in \(\mathcal{X}  \) is defined by
\[ 
\mathbf{A}_{\infty} \boldsymbol{\rho} = \{- \boldsymbol{\Delta} \mathbf{w} \mid \mathbf{w}_i \in D(L), \boldsymbol{\rho}_i = \text{sign}(\mathbf{w}_i) \text{ a.e. on } \Omega, \quad i = 0, 1, \dots, K \}.
\]
Based on \cite{brezis1973semi}, \(A_m\) is m-accretive in \(\Omega\) for \(m \in [1, \infty]\), then
\[ 
(I + A_m)^{-1} \rho \to (I + A_\infty)^{-1} \rho \quad \text{in } \Omega \quad \text{as } m \to \infty.
\]
Using the orthogonality of basis functions, thanks to \cite{benilan1988some}, 
one can show that the vectorized operator satisfies
\begin{equation}\label{eq:invConv}
	(I + \mathbf{A}_m)^{-1} \boldsymbol{\rho} \to (I + \mathbf{A}_\infty)^{-1}  \boldsymbol{\rho} \quad \text{in } \mathcal{X} \quad \text{as } m \to \infty.
\end{equation}
Furthermore, as for the semigroups generated by these operators \cite{benilan1989limit}, we have
\begin{equation}\label{eq:semigroup}
	e^{-t\mathbf{A}_m} \circ \mathbf{f} \to e^{-t\mathbf{A}_\infty} \circ \mathbf{f} \quad \text{in } \mathcal{X} \quad \forall\, t>0,
\end{equation}
where \( \mathbf{f} \) is the coefficient vector corresponding to the initial datum \( f \). Apply \cite{benilan1996singular} Theorem 1, we can prove Theorem 2.

\end{document}